\documentclass[twoside]{amsart}
\usepackage{myamsart}
\input{mydef}
\providecommand{\eqnright}[1]{\qquad\qquad  #1}

\providecommand{\Ker}{\object{Ker}}
\providecommand{\Ess}{\object{Ess}}
\providecommand{\CH}{\object{CH}}
\providecommand{\bos}{{\mathcal B}}

\begin{document}
\title[Relative Convolutions. I]{Relative Convolutions. I\\
\itshape\Large Properties and Applications}

\author{Vladimir V. Kisil}
\email{kisilv@maths.leeds.ac.uk}
\address{%
School of Mathematics,
University of Leeds,
Leeds LS2\,9JT,
UK}

\thanks{This work was partially supported
by CONACYT
Project 1821-E9211, Mexico.}
\thanks{On leave from the Odessa State
University.}
\date{September 2, 1994}
\begin{abstract}
To study operator algebras with symmetries in a wide sense
we introduce a notion of {\em relative convolution
operators\/} induced by a Lie algebra. Relative convolutions recover
many
important classes of operators, which have been already studied
(operators
of multiplication, usual group convolutions, two-sided convolution etc.)
and their different combinations. Basic properties of relative
convolutions are given and a connection with usual convolutions is
established.

Presented examples show that relative convolutions provide us
with a base for systematical applications of harmonic analysis
to PDO theory, complex and hypercomplex analysis, coherent states,
wavelet transform and quantum theory.

\keywords{Lie groups and algebras, convolution operator, representation
theory, Heisenberg group, integral representations,
Hardy space, Szeg\"o projector, Toeplitz
operators, Fock space, Segal--Bargmann
space, Bargmann projector, Dirac equation, Clifford analysis, coherent
states, wavelet transform, quantization.}
\AMSMSC{45P05}{43A80, 22E45, 22E60, 32M05, 81R30, 81S99}
\end{abstract}
\maketitle
\newpage
\tableofcontents
\section{Introduction}
Convolution operators and different operators associated with them are
very important in mathematics (for example, in complex
analysis~\cite{GreSte77}) and
physics (for example, quantum mechanics~\cite{Folland89}
and~\cite[Chap.~3]{Woodhouse80}).
The
fundamental role of such operators may be easily explained. Indeed,
a notion of symmetries and group transforms invariance are the basis of
the contemporary science. It is well known, that operators, which are
invariant under a transitive group operation, may be realized as
convolution
operators on the group.

However, there are some limitations for an
application of convolution operators:
\begin{itemize}
\item It is not a rare case when a operators symmetry group and
functions
domain have different dimensions.
\item  Group convolution operators are generated by
transformations of function domains. However, mathematical objects
are often connected not only with domain transforms but also with alterations
of function range or even both of them.
\end{itemize}
The paper introduces a notion of
{\em relative convolutions\/}, which allows us to overcome these
limitations. Naturally,  such an interesting object cannot
be totally unknown in mathematics and we will give a short description
of connected ones in Remark~\ref{re:origin}. In this paper we work
only with continuous groups of symmetry. Discrete groups or their
mixtures with continuous ones may be also considered.

A geometric group structure and related objects (especially the
non-commutative
Fourier transform) are obviously determining for properties of
convolution algebra. If the group structure is used not
only
as a basis for particular calculations but also in a more general
framework, then
it is possible to find the proper level of generality for obtained
results\footnote{See, for example, the deduction of a PDO--form for
convolutions
on step 2 nilpotent Lie groups in~\cite[p.~9--10]{MTaylor84} or
Theorem~\ref{th:projector}.}.

Be found, that geometry of Lie groups and
algebras jointly with properties of kernels determine representations of
relative convolution algebras. Namely, all representations of relative
convolution algebras will be induced by selected representations of
corresponding
Lie algebras (Lie groups) and the selection of representations  will
depend on kernels properties. The representation theory of Lie algebras
is complicated and still unsolved completely. But, due to the author
personal interest, the main examples are nilpotent Lie algebras fully
described by the Kirillov theory (see~\cite{Kirillov62}
or~\cite[Chap.~6]{MTaylor86}). Considered examples will show that such a
restriction still leaves enough space for very interesting applications.

 The layout is as follows.

In Section~\ref{se:relative} we give necessary notations, introduce the
main object of the paper --- {\em relative convolution\/} --- and show
basic examples. Relative convolutions recover many
important classes of operators, which have been already studied
(operators of
multiplication, usual group convolutions, two-sided convolutions etc.)
and their different combinations.

Basic properties of relative convolutions are given in
Section~\ref{se:properties}. We prove a formula for the composition of
relative convolutions and deduce from it that {\em an algebra of
relative convolutions induced by a Lie algebra \algebra{g}
is a representation of the algebra of group convolutions on the Lie
group $\object{Exp}\algebra{g}$}. We also state the universal role PDO
and the Heisenberg group in the  theory of relative
convolution operators.

Relative convolutions provide us with a tool for the systematic usage of
harmonic analysis in different fields of pure and applied mathematics.
This statement will be
illustrated in Sections~\ref{se:analysis} and~\ref{se:physics} by many
examples, which
show applications of relative
convolutions to the theory of PDO, complex and hypercomplex analysis,
coherent states, wavelets and
quantum mechanics.

In Section~\ref{se:coherent} we will describe our view on a notion of
coherent states. The main observation is: a direct employment of group
structures of coherent states gives us the uniform theory for
the Bargmann, Bergman and Szeg\"o projectors at the
Segal--Bargmann (Fock), Bergman and Hardy spaces respectively.
Applications to wavelets (and other) theory are also possible.

The author gratefully acknowledge his inspiration by
papers~\cite{Dynin75,Dynin76,Howe80a,Howe80b,MTaylor84}.
Besides,  there is some interference with the coherent recent paper of
Folland~\cite{Folland94}. The author is also grateful to
Dr.~V.~V.~Kravchenko and
Prof.~N.~L.~Vasilevski  for helpful discussions.

\section{Relative Convolutions}\label{se:relative}
\subsection{Definitions and Notations}\label{ss:notation}
Let $G$ be a connected simply connected Lie group, let $\algebra{g}$ be
its finite-dimensional Lie algebra. By the usual way we can identify
\algebra{g} with \Space{R}{N} for $N=\dim\algebra{g}$ as vector spaces.
The exponential map $\exp: \algebra{g}\rightarrow
G$~\cite[\S~6.4]{Kirillov76} identifies the group $G$ and its algebra
\algebra{g}. So
$G$  as a real $C^\infty$-manyfold has the dimension $N$. Let fix a
frame $\{X_j\}_{1\leq j\leq N}$ of \algebra{g}. Via the exponential map we
can write a decomposition $g=\sum_1^N x_j X_j$ in {\em exponential
coordinates\/} for every $g\in G$.
The group law on $G$ in the exponential coordinates may be expressed via
the {\em Campbell--Hausdorff formula\/}~\cite[\S~6.4]{Kirillov76}:
\begin{equation}\label{eq:Campbell-Hausdorff}
g*h=\sum_{m=1}^\infty\frac{(-1)^{m-1}}{m}
\sum_{\vbox{\hbox{\hfil\makebox[\width]{\ $\scriptstyle k_j+l_j\geq
1$}\hfil}
\hbox{\makebox[\width]{$\scriptstyle k_j\geq 1,\, l_j\geq 1$}}}}
\frac{[x^{k_1}y^{l_1}\ldots x^{k_m}y^{l_m}]}{k_1!l_1!\ldots k_m!l_m!},
\end{equation}
where $[x_1x_2\ldots x_n]=\frac{1}{n}[\ldots[[x_1,x_2],x_3],\ldots,x_n]$
and
\begin{displaymath}
g=\sum_1^N x_j X_j,\ h=\sum_1^N y_j X_j.
\end{displaymath}
It seems reasonable to
introduce a short notation for the right side of
{}~\eqref{eq:Campbell-Hausdorff}, we select $\CH[\sum_1^N x_j X_j,\sum_1^N
y_j X_j]$ or $\CH[ x,y]$ (if the frame is obvious).

Let $S$ be a set and let be defined an operation $G: S\rightarrow S$ of
$G$ on  $S$. If we fix a point $s\in S$ then the set of elements
$G_s=\{g\in G\such g(s)=s\}$ obviously forms the {\em isotropy
(sub)group of
$s$
in $G$\/}~\cite[\S~I.5]{Lang69}. There is an equivalence relation on $S$,
say, $s_1\sim
s_2 \Leftrightarrow \exists g\in G: gs_1=s_2$, with respect to which $S$
is a disjoint union of distinct orbits~\cite[\S~I.5]{Lang69}. Thus from now on,
without lost of a generality, we assume that the
operation of $G$ on $S$ is {\em transitive\/}, i.~e. for every $s\in S$
we have
\begin{displaymath}
Gs:=\relstack{\bigcup}{g\in G} g(s)=S.
\end{displaymath}
If $G$ is a Lie group then the {\em homogeneous space\/} $G/G_s$ is also a Lie
group for every $s\in
S$. Therefore the one-to-one mapping $G/G_s \rightarrow S: g\mapsto g(s)$
induces a
structure of $C^\infty$-manifold on $S$. Thus the class
$\FSpace[\infty]{C}{0}(S)$ of smooth functions with compact supports on
$S$ has the evident definition.

A smooth measure $d\nu$ on $S$ is called  {\em invariant\/} (the {\em Haar
measure\/}) with respect to an operation of $G$ on $S$ if
\begin{equation}
\int_S f(s) d\nu(s)\nonumber= \int_S f(g(s)) d\nu(s),
\mbox{
for all } g\in G,\ f(s)\in\FSpace[\infty]{C}{0}(S).\label{eq:invar-m}
\end{equation}
The Haar measure always exists and is uniquely defined up to a scalar
multiplier~\cite[\S~0.2]{MTaylor86}. An equivalent formulation
of~\eqref{eq:invar-m} is: {\em $G$ operate on $\FSpace{L}{2}(S,d\nu)$ by
unitary operators\/}. We will transfer the Haar measure
$d\mu$ from $G$ to \algebra{g} via the exponential map $\exp:
\algebra{g}\rightarrow G$ and will call it as the {\em invariant measure
on a Lie algebra \algebra{g}\/}.

Now we can define an operation of a Lie algebra \algebra{g} in the space
$\FSpace[\infty]{C}{0}(S)$ induced by an operation of $G$ on $S$. Let
$X\in
\algebra{g}$ and $f(s)\in
\FSpace[\infty]{C}{0}(S)$ then
\begin{equation}\label{eq:alg-limit}
\lim_{t\rightarrow 0}\frac{f(e^{tX}s)-f(s)}{it} \in
\FSpace[\infty]{C}{0}(S),
\end{equation}
where $\exp: X\mapsto e^{tX}$ is the exponential map $\exp:
\algebra{g}\rightarrow
G$. The value of limit~\eqref{eq:alg-limit} will be denoted by
$[Xf](s)$.
If $S$ is equipped by a measure $d\nu$  we can define an {\em adjoint\/}
operation $X^*$ of \algebra{g} on
$\FSpace{L}{2}(S,d\nu)$
by natural formula $\scalar{f(s)}{[X^* g](s)}:=\scalar{[Xf](s)}{g(s)}$.
The invariance~\eqref{eq:invar-m} of $d\nu$ may be reformulated in the
terms of the Lie algebra \algebra{g} as follows:
\begin{displaymath}
X_j^* =X_j,
\end{displaymath}
for every $X_j,\ 1\leq j\leq N$ from a frame of \algebra{g}.
Summarizing, $X: f\rightarrow [Xf]$ is a selfadjoint (possibly
unbounded) operator in $\FSpace{L}{2}(S,d\nu)$ with invariant measure
$d\nu$.

In general, we will
speak on an {\em operation\/} of a Lie algebra \algebra{g} on a manifold
$S$ with a measure $d\nu$ if there is a linear representation of
\algebra{g} by selfadjoint operators on  the linear
space $\FSpace{L}{2}(S,d\nu)$. As usual, if $X$ is selfadjoint on
$\FSpace{L}{2}$ then $e^{itX}$ is unitary on \FSpace{L}{2}. Clearly,
every
operation of a Lie group on $S$ induces an operation of the
corresponding Lie
algebra, but inverse is not true generally speaking (see
Example~\ref{ex:multiplication}). Therefore we will based on the notion
of a {\em Lie algebra\/} operation.

The succeeding object will be useful in our study of convolution
algebras
representations. We define the {\em kernel\/} of operation of
\algebra{g} on $S$ as follows:
\begin{equation}\label{eq:def-kernel}
\Ker(\algebra{g},S)=\{X\in \algebra{g} \such [Xf](x)=0 \mbox{ for all }
f(x)\in\FSpace[\infty]{C}{0}(S)\}.
\end{equation}
If an operation of
\algebra{g} on $S$ is induced by an operation of $G$ on $S$ then
$X\in\Ker(\algebra{g},S)$ if and only if for $e^X\in G$ and for any
$s\in S$ we have $e^X(s)=s $.
There is no doubt that $\Ker(\algebra{g},S)$ is a two-sided ideal of
\algebra{g} (and therefore a linear subspace). Thus we can introduce the
quotient Lie algebra
$\Ess(\algebra{g},S)=\algebra{g}/\Ker(\algebra{g},S)$
(see~\cite[\S~6.2]{Kirillov76}). An induced action of
$\Ess(\algebra{g},S)$ on $S$
is evidently specified.

Now we can describe the main object of the paper.
\begin{defn}\label{de:relative}
Let we have a (selfadjoint) operation of a Lie algebra \algebra{g} on
$S$ (possibly
induced by an operation of a Lie group $G$ on a set $S$) and let
$\{X_j\}_{1\leq j\leq N}$ be a fixed frame of \algebra{g}. We will
define
the {\em operator of relative convolution\/} $K$ induced by \algebra{g}
on
$E=\FSpace[\infty]{C}{0}(S)$ with a kernel $k(x)\in
\FSpace[\infty]{C}{0}(\Space{R}{N})$ by the formula
\begin{equation}\label{eq:relative}
K=(2\pi)^{-N/2}
\int_{\Space{R}{N}} \widehat{k}(x_1,x_2,\ldots,x_N)\,e^{i\sum_1^N x_j X_j}\,
dx,
\end{equation}
where integration is made with respect to an invariant measure on
$\algebra{g}\cong\Space{R}{N}$.
Here $\widehat{k}(s)$ is the Fourier transform of the function $k(s)$
over \Space{R}{N}:
\begin{displaymath}
\widehat{k}(x)= (2\pi)^{-N/2}\int_{\Space{R}{N}}k(y)\,e^{-iyx}\,dy.
\end{displaymath}
\end{defn}
\begin{rem} \label{re:origin}
This definition has an origin at the {\em Weyl functional
calculus\/}~\cite{HWeyl}
(or the Weyl quantization procedure), see Example~\ref{ex:pdo} for
details. Feynman in~\cite{Feynman51} proposed its extension --- the {\em
functional calculus of ordered operators\/} --- in a very similar way.
Anderson~\cite{Anderson69} introduced a generalization of the Weyl
calculus for arbitrary set of self-adjoint operators in a Banach space
exactly by formula~\eqref{eq:relative}. A description of different
operator calculuses may be found in~\cite{Maslov73}. But it was shown by
R.~Howe~\cite{Howe80a,Howe80b}, that success of the original Weyl
calculus is intimately connected with the structure of the Heisenberg
group and its different representations. Thus one can obtain a
new fruitful branch in this direction making an assumption, that
the operators $X_j$
in~\eqref{eq:relative} are not arbitrary but are connected with some
group structure. Such a treatment for the Heisenberg group and multipliers
may be found at the Dynin paper~\cite{Dynin75}. M.~E.~Taylor
in~\cite{MTaylor84} introduced ``smooth families of convolution
operators'', which technically coincide with relative convolutions in
many important cases. However, full symmetry groups of such smooth
families are not clear and thus cannot be exploited. Efforts to
study simultaneously both the group action and operators of multiplication
bring also a very similar object at the recent
paper\footnote{In the mentioned paper consideration is restricted to
nilpotent step 3 Lie algebras.} of Folland~\cite{Folland94}. But general
consideration of relative convolutions  seems to be new (as well as
systematical applications to new and already solved problems, see
Sections~\ref{se:analysis} and~\ref{se:physics}).
\end{rem}
Due to paper~\cite{Anderson69} our definition is correct for any
representation of a Lie algebra in a Banach space. Thus we can use
similar
\begin{defn}\label{de:relative1}
Let we have a selfadjoint representation of a Lie algebra \algebra{g} on
a Banach space $B$ and let
$\{X_j\}_{1\leq j\leq N}$ be operators represented a fixed frame of
\algebra{g}. We will define the {\em operator of relative convolution\/}
$K$ induced by \algebra{g} on $B$ by the formula
\begin{equation}\label{eq:relative1}
K=(2\pi)^{-N/2}
\int_{\Space{R}{N}} \widehat{k}(x_1,x_2,\ldots,x_N)\,e^{i\sum_1^N x_j X_j}\,
dx.
\end{equation}
Here integration is made with respect to an invariant measure on
$\algebra{g}\cong\Space{R}{N}$.
\end{defn}
\begin{rem}\label{re:function}
As was shown at~\cite{Anderson69}, formula~\eqref{eq:relative1} may be treated
as a definition of a function $k(X_1,X_2,\ldots,X_N)$ of operators $X_j$. See
Remark~\ref{re:riesz} for an alternative Riesz-like definition of
relative convolutions.
\end{rem}
We have
defined relative convolutions only for a very restricted class of
kernels $k(x)$ and a specific space $E$. Of course, both of them may be
enlarged in
the proper context. Another interesting modification may be required by
a consideration of discrete groups and their actions. It may be achieved
by a replacing integration in~\eqref{eq:relative} by summation or
integration by discrete measures. But in this paper we will not
consider such cases.

As we
will see, relative convolutions are naturally defined not only for group
operation on {\em functions domains\/} but also on {\em ranges of
function\/} (see Examples~\ref{ex:multiplication}, \ref{ex:fock}
and~\ref{ex:clifford}).

\subsection{Basic Examples}
The
following Examples make clear the relationships between the relative
convolutions and usual ones.
\begin{example}\label{ex:euclid-conv}
Let $G=\Space{R}{N}$ and $G$ operates on $S=\Space{R}{N}$ as a group of
Euclidean shifts $g: y\rightarrow y+g$. The algebra \algebra{g} consists
of selfadjoint differential operators spanned on the frame
$ X_j^e= \frac{1}{i} \frac{\partial}{\partial y_j},\ 1\leq j\leq N$.
Then:
\begin{eqnarray*}
[Kf](y)&=& (2\pi)^{-N/2}\int_{\Space{R}{N}} \widehat{k}(x)\,e^{i\sum_1^N
            -x_j (\frac{1}{i}\frac{\partial }{\partial y_j})}f(y)\,dx\\
       &=& (2\pi)^{-N/2}\int_{\Space{R}{N}} \widehat{k}(x)\,e^{-\sum_1^N
            x_j\frac{\partial }{\partial y_j}}\,f(y)\,dx\\
       &=&(2\pi)^{-N/2} \int_{\Space{R}{N}} \widehat{k}(x)\,f(y-x)\,dx .
\end{eqnarray*}
Otherwise, an operator of relative convolution with  a kernel $k(x)$
evidently coincides with a usual (Euclidean) convolution on
\Space{R}{N} with the kernel $\widehat{k}(x)$.
\end{example}
This Example can be obviously generalized for an arbitrary Lie group $G$
and
the set $S=G$ with the natural left\footnote{For a commutative
group G the
left and the right operations are the same.}  (or right) operation of
$G$ on
$S=G$
\begin{displaymath}
G: S \rightarrow GS\ (SG): g' \mapsto g^{-1}g'\ (g' \mapsto g'g),\ g\in G,\
g'\in S=G.
\end{displaymath}
It is clear, that relative convolutions will coincide with the left
(right) group
convolutions on $G$ (with, may be, transformed kernels).
\begin{example}\label{ex:multiplication}
Let $G=\Space{R}{N}$ with the Lie algebra spanned on operators of
multiplication $X_j=M_{y_j}$ by $y_j,\ 1\leq j \leq N$ on a space of
functions on
$S=\Space{R}{N}$. Then a relative convolution $K$ with a kernel $k(x)$
\begin{eqnarray*}
[Kf](y)&=&(2\pi)^{-N/2}\int_{\Space{R}{N}} \widehat{k}(x)\,e^{ixy}\,f(y)\,dx\\
       &=& k(y)f(y)
\end{eqnarray*}
is simply an operator of multiplication by $k(y)$. In this case the
operation of the Lie algebra \algebra{g} on $S$ is not generated by an
operation of a Lie group on $S$. Of course, it is possible to
establish a connection with Example~\ref{ex:euclid-conv} through the
Fourier transform, but it is not so simple in other cases. Particularly
it will happen when one generalizes this Example by use more complicated
transformations of range than a multiplication by scalars (see
Example~\ref{ex:clifford}) or simultaneously applies transformations from this
and
the previous Examples (see Example~\ref{ex:pdo}).
\end{example}
\begin{example}\label{ex:two-sided}
Let \Heisen{n} be the Heisenberg
group~\cite{Dynin75,MTaylor84,MTaylor86}.
 The Heisenberg group \Heisen{n} is a step 2 nilpotent Lie group. As a
$C^\infty -$manifold it coincides with $\Space{ R}{2n+1}$. If an
element of it is given in the form $g=(u,v)\in \Space{H }n$, where
 $u\in \Space{ R}{}$ and $v=(v_1,\dots, v_n)\in \Space{ C}n$, then the
group law on $
\Space{ H}{n}$ can be written as
\begin{equation}\label{eq:h-low}
(u,v)*(u',v')=\left(u+u'-\frac{1}{2} \object{Im} \sum_1^n v'_k
\bar{v}_k,\
v_1+v'_1,\ldots, v_n+v'_n\right) .
\end{equation}
We single out on $\Space{ H}n$ the group of nonisotropic dilations
 $\{\delta_\tau\}$, $\tau\in \Space{ R}{}_+$:
\begin{displaymath}
\delta_\tau(u,v)=(\tau^2u,\tau v).
\end{displaymath}
Functions with the property
\begin{equation}\label{eq:homohen}
(f\circ \delta_\tau)(g)=\tau^kf(g)
\end{equation}
will be called {\em $\delta_\tau$-homogeneous functions of degree
$k$\/}.

The left and right Haar measure on the Heisenberg group coincides with
the Lebesgue
measure.  Let us introduce the  right
$\pi_{r}$ and the left $\pi_{l}$ regular representations of \Heisen{n}
on
$\FSpace{L}{2}(\Heisen{n})$:
 \begin{eqnarray}
{} [\pi_{l} (g) f] (h) &=& f(g^{-1}*h),\label{eq:lf-shift}\\
{} [\pi_{r} (g) f] (h) &=& f(h*g). \label{eq:rt-shift}
 \end{eqnarray}
Thereafter,
$G=\Heisen{n}\times\Heisen{n}$, $S=\Heisen{n}$ and an operation of $G$
on
$S$ is defined by the {\em two-sided\/} shift
\begin{displaymath}
G: s \rightarrow
g^{-1}_1*s*g_2,\ s\in S=\Heisen{n}, (g_1,g_2)\in G=\Heisen{n}
\times\Heisen{n}.
\end{displaymath}
The Lie algebra of $G$ is the direct sum
$\algebra{g}=\algebra{h}_n\oplus\algebra{h}_n$ of two copies of the Lie
algebra $\algebra{h}_n$ of the Heisenberg group \Heisen{n}. If in
$S\cong\Space{R}{N}$ we define Cartesian coordinates $y_j,\ 0\leq j<
N=2n+1$, then a frame of \algebra{g} may be written as follows
\begin{eqnarray}
\displaystyle X^l_0=\frac{1}{i}\frac{\partial }{\partial y_0}, &
\displaystyle X^l_j=\frac{1}{i}\left(\frac{\partial}{\partial
y_{j+n}}-2y_j\frac{\partial }{\partial y_0 }\right), &\displaystyle
X^l_{j+n}=\frac{1}{i}\left(\frac{\partial}{\partial
y_{j}}+2y_{j+n}\frac{\partial
}{\partial y_0 }\right),\nonumber \\
  \label{eq:frame-left}\\
\displaystyle X^r_0=-\frac{1}{i}\frac{\partial }{\partial
y_0},&\displaystyle  X^r_j=\frac{1}{i}\left(\frac{\partial}{\partial
y_{j+n}}+2y_j\frac{\partial }{\partial y_0 }\right), &\displaystyle
X^r_{j+n}=\frac{1}{i}\left(\frac{\partial}{\partial y_{j}}-
2y_{j+n}\frac{\partial
}{\partial y_0 }\right),\nonumber \\
\label{eq:frame-right}
\end{eqnarray}
where $1\leq j\leq n$. Vector fields $X^l_j,\ 0\leq j < N$ generate the
left shifts~\eqref{eq:lf-shift} on $S=\Heisen{n}$ and $X^r_j,\ 0\leq j <
N$ generate the
right ones~\eqref{eq:rt-shift}. $X^l_j$
(correspondingly $X^r_j)$  satisfy to the famous Heisenberg commutation
relation
\begin{equation}\label{eq:heisen-comm}
[X^{l(r)}_j,X^{l(r)}_{j+n}]=-[X^{l(r)}_{j+n},X^{l(r)}_{j}]= X^{l(r)}_{0}
\end{equation}
and all other commutators being zero. Particulary, all $X^l_j$ commute
with all $X^r_j$.

Vector fields $X^l_0$ and $-X^r_0$ are different as elements of
\algebra{g}, but have coinciding operations on functions. It is easy to
see, that $\Ker(\algebra{g},S)$ is a linear span of the vector
$X^l_0+X^r_0$ and $\Ess(\algebra{g},S)=\Heisen{2n}$.
Now formula~\eqref{eq:relative} with kernel $k(x),\
x\in\Space{R}{4n+2} $ defines the {\em two-sided convolutions\/} studied
in~\cite{VasTru94}:
\begin{eqnarray*}
[Kf](y)&=&(2\pi)^{-N}\int_{\Space{R}{4n+2}}
\widehat{k}(x)\,e^{i\sum_1^{2N}x_j X_j}\,
f(y)\,dx\\
       &\stackrel{(*)}{=}&(2\pi)^{-N}\int_{\Space{R}{4n+2}}
\widehat{k}(x)\,e^{i\sum_0^{N-1}x'_j X^l_j}\,
e^{i\sum_0^{N-1}x''_j X^r_j}\, f(y)\,dx\\
       &=&(2\pi)^{-N}\int_{\Heisen{n}} \int_{\Heisen{n}} \widehat{k}(x',x'')\,
           f(x'^{-1}*y*x'')\,dx'dx''\\
       &=&(2\pi)^{-N} \int_{\Heisen{n}} \int_{\Heisen{n}}
\widehat{k}(x',x'')\pi_l(x')\pi_r(x'')\,dx'dx''\,f(y).
\end{eqnarray*}
Transformation $(*)$ is possible due the commutativity of $X_j^l$ and
$X_j^r$.

Reduction of two-sided convolutions to the usual group ones was done
in~\cite{Kisil94f} in a way very similar to the present consideration
(see Corollary~\ref{co:two-sided}).
\end{example}

The Heisenberg group here may be substituted by any non-commutative
group
$G$ (for a commutative group the left and the right shifts are the same)
and we will obtain two-sided convolutions on $G$.

These basic Examples form a frame for other ones: the number of
examples may be increased both by simple compositions of
convolutions from Examples~\ref{ex:euclid-conv}--\ref{ex:two-sided} and
by alterations of considered groups.

\section{Basic Properties}\label{se:properties}
The main purpose of the present introduction to relative convolutions is
a sharp extension of harmonic analysis applications. Through relative
convolutions we can transfer the knowledge on Lie groups and their
representations to the different operator algebras. Properties of
relative convolutions themselves are strictly depending on a structure of
a concrete group. Thus, it seems unlikely that one can say too much
about them in general. Nevertheless, results collected in this Section
establish important properties of relative convolutions and will be
forceful in the future.

$\Ker(\algebra{g},S)$ operates on $S$ trivially therefore we are
interested to understand the effective part of a relative convolution
operator. Let $m:=\dim\Ker(\algebra{g},S)>0$. We have a decomposition
$\algebra{g}=\Ker(\algebra{g},S)\oplus\Ess(\algebra{g},S)$ as linear
spaces.
Let us introduce a frame $\{X_j\}_{1\leq j \leq N}$ of \algebra{g} such
that $X_j \in \Ker(\algebra{g},S), 1\leq j \leq
m$ and we assume that
operator~\eqref{eq:relative} is written through such a frame (it clearly
can be obtained by a change of variables). Then the following Lemma is
evident.
\begin{lem}\label{le:decomposition}\textup{(Effective decomposition)}
Operator~\eqref{eq:relative} of relative convolution is equal to
an operator of relative convolution generated by the induced operation
of $\Ess(\algebra{g},S)$ on $S$ with the kernel $k'(x'')$ defined by the
formula
\begin{equation}\label{eq:decomposition}
\widehat{k'}(x'')=(2\pi)^{-m/2}
\int_{\Space{R}{m}}\widehat{k}(x',x'')\,dx',
\end{equation}
where $x'\in\Space{R}{m},\ x''\in\Space{R}{N-m}$.
\end{lem}
The effective decomposition allows us to consider always a basic
case of $\Ker(\algebra{g},S)=0,\ \Ess(\algebra{g},S)=\algebra{g}$. After
that the general case may be obtained by the obvious modification.

Operators of relative convolution clearly form an algebra and we will
denote by $(k_2*k_1)(x)$ the kernel of composed operator of two relative
convolutions with
kernels $k_1(x)$ and $k_2(x)$. Of course, $k_2*k_1\neq
k_1*k_2$ generally speaking.
\begin{prop}\textup{(Composition formula)} \label{pr:composition} We
have
\begin{equation}\label{eq:composition}
(k_2*k_1)(x)=(2\pi)^{-N/2}\int_{\Space{R}{N}}k_2(y)\,k_1(\CH[-y,x])\,dy,
\end{equation}
where $\CH[-y,x]=\CH[y^{-1},x]$ is given by~\eqref{eq:Campbell-Hausdorff}.
\end{prop}
\begin{proof} For kernels $k_1, k_2$ satisfying the Fubini theorem we
can change theintegration order and obtain
formula~\eqref{eq:composition}:
\begin{eqnarray*}
K_2 K_1&=& (2\pi)^{-N/2}\int_{\Space{R}{N}}k_2(y)\,e^{i\sum_1^N y_j X_j}\,dy
         \times (2\pi)^{-N/2}\int_{\Space{R}{N}}k_1(z)\,e^{i\sum_1^N z_j
X_j}\,dz\\
 &=& (2\pi)^{-N}\int_{\Space{R}{N}}\int_{\Space{R}{N}}k_2(y)\,k_1(z)\,
    e^{i\sum_1^N y_j X_j}e^{i\sum_1^N z_j X_j}\,dy \,dz\\
 &=& (2\pi)^{-N}\int_{\Space{R}{N}}\int_{\Space{R}{N}}
k_2(y)\,k_1(z)\,e^{i\CH[\sum_1^N y_j X_j,\sum_1^N z_j X_j]}\,dy \,dz\\
 &=& (2\pi)^{-N/2}\int_{\Space{R}{N}}(2\pi)^{-N/2}\int_{\Space{R}{N}}
k_2(y)\,k_1(\CH[-\sum_1^N y_j X_j,\sum_1^N x_j X_j])\,dy\\
 &&\eqnright{\times  e^{i\sum_1^N x_j X_j}\,dx,}
\end{eqnarray*}
where $e^{i \sum_1^N x_j X_j}=e^{i \sum_1^N y_j X_j}e^{i \sum_1^N z_j
X_j}$. Making the change of variables we used the invariance of the measure on
\algebra{g}.
\end{proof}
Note, that for an operation \algebra{g} on $S=G$ induced by
the operation of $G$ on $S$ formula~\eqref{eq:composition} gives us the
usual {\em group convolution\/} on $S=G$. Thus we can speak
about a {\em convolution-like\/} calculus of relative convolutions. This
remark gives us the following important result describing the nature of
relative convolution algebras\footnote{Relative convolutions obviously
include usual group convolutions and, on the other hand, due to
Theorem~\ref{th:nature}
they may be treated just as representations of group convolutions. Thus
they are not {\em generalized\/} convolutions but simply {\em
relative\/} convolutions.}.
\begin{thm}\label{th:nature}
Let \algebra{g} be an algebra Lie operating on a set $S$ and let $G$ be
the exponential Lie group
of $\Ess(\algebra{g},S)$. Let also \algebra{G} be an algebra of
relative convolutions induced
by \algebra{g} on $S$ and let $\widehat{\algebra{G}}$ be a group
convolution algebra on $G$. Then \algebra{G} is a linear representation
of $\widehat{\algebra{G}}$.
\end{thm}
\begin{rem}Due to Theorem~\ref{th:nature} we can give an alternative
definition of a relative convolution algebra, namely, {\em an operator
algebra \algebra{G} is a relative convolution algebra induced by a Lie
algebra \algebra{g}, if all representations of \algebra{G} are
subrepresentations of the convolution algebra on group
$\object{Exp}\algebra{g}$\/}.
\end{rem}
 Theorem~\ref{th:nature} is the main tool for
applications of harmonic analysis to every problem where relative
convolutions occur.
Particulary, it gives us an easy access to many conclusions obtained
by direct calculations. Now we illustrate this by applications to
two-sided convolutions from Example~\ref{ex:two-sided}.
\begin{cor}\upshape{\cite{Kisil94f}}\label{co:two-sided} An algebra of
two-sided convolutions on the
Heisenberg group \Heisen{n} is a representation of the group convolution
algebra on \Heisen{2n}.
\end{cor}
Explicit descriptions of the established representation depend from
properties of relative convolutions kernels. In the next paper will be
shown that for kernels from
$\FSpace{L}{1}(\Space{R}{N}\times\Space{R}{N})$ there
is a one-to-one correspondence between representations of two-sided
convolutions on \Space{R}{N}  and usual one-sided convolution on the
$\object{Exp}\Ess(\algebra{g}\times\algebra{g},\Space{R}{N})$ (for the
Heisenberg group it was calculated in~\cite{Kalyuzhny93}). If kernels
have a symmetry group then some representations of
$\object{Exp}\Ess(\algebra{g}\times\algebra{g},\Space{R}{N})$ vanish,
see, for example case of two-sided convolutions on \Heisen{n} with
homogeneous kernels~\eqref{eq:homohen} in~\cite{VasTru94}.

There is another simple but important corollary of
Theorem~\ref{th:nature}
\begin{cor}\label{co:auto}
Let $\psi: \algebra{g} \rightarrow\algebra{g}$ be an automorphism of
\algebra{g} as a Lie algebra. Let
$\Psi:\algebra{G}\rightarrow\algebra{G}$ be
a mapping defined by the rule $\Psi: K \mapsto \Phi(K)$, where $K$ is a
relative convolution with a kernel $k(x)$, $\Psi(K)$ is a relative
convolution with the kernel $k(\psi x)\, J^{1/2}{\psi}(x)$ and
$J{\psi}(x)$ is the Jacobian of $\psi$ at the point $x$. Then $\Psi$ is
an automorphism of algebra \algebra{G}.
\end{cor}
\begin{proof} The key point of the proof is the invariance
of~\eqref{eq:composition} under $\Psi$ (the rest is almost
evident). This invariance follows from a simple change of variables in
the integral~\eqref{eq:relative}:
\begin{eqnarray*}
\Psi[k_1*k_2](x)&=&(2\pi)^{-N/2}\int_{\Space{R}{N}}k_2(y)\,k_1(\CH[-
y,\psi(x)])\,   J^{1/2}{\psi(x)}\,dy\\
       &=&(2\pi)^{-N/2}\int_{\Space{R}{N}}k_2(y)\,k_1(\CH[-\psi(\psi^{-1}
(y)),\psi(x)])\,   J^{1/2}{\psi(x)}\,dy\\
       &=&(2\pi)^{-N/2}\int_{\Space{R}{N}}k_2(\psi(y'))\,k_1(\CH[-
\psi(y'),\psi(x)])\,   J^{1/2}{\psi(x)}\,dy\\
&\stackrel{(*)}{=}&(2\pi)^{-N/2}\int_{\Space{R}{N}}k_2(\psi(y'))\,k_1(\psi(\CH[-
y',x])) \, J^{1/2}{\psi(x)}\,dy\\
       &=&(2\pi)^{-N/2}\int_{\Space{R}{N}}k_2(\psi(y')) \, J^{1/2}{\psi(y')}\,
k_1(\psi(\CH[-y',x]))\\
       &&\eqnright\times  J^{1/2}{\psi(y'*x)} \,  J^{-1}{\psi}(y')\,dy\\
       &=&(2\pi)^{-N/2}\int_{\Space{R}{N}}\, \Psi k_2(y')\,\Psi k_1(\CH[-
y',x]) \,dy'\\
       &=& [\Psi k_2*\Psi k_1](x),
\end{eqnarray*}
where $y'=\psi^{-1}(y),\ y=\psi(y')$. Here transform $(*)$ is possible due to
the automorphism
property of $\psi$. We also employ the identity
\begin{displaymath}
J^{1/2}{\psi(y')}\,  J^{1/2}{\psi(y'*x)}\, J^{-1}{\psi}(y')=
J^{1/2}{\psi(x)},
\end{displaymath}
which follows from the chain rule.
\end{proof}
Next Lemma evidently follows from Definition~\ref{de:relative1} and
Theorem~\ref{th:nature}.
\begin{lem}\label{le:category}
Let an algebra $\widetilde{\algebra{G}}$ be a representation of an
algebra \algebra{G} of relative convolutions induced by a Lie algebra
\algebra{g}. Then $\widetilde{\algebra{G}}$ is also the algebra of
relative convolutions induced by \algebra{g}.

Otherwise, algebras of relative convolutions induced by a Lie algebra
\algebra{g} form a closed sub-category~$\mathcal G$ of the category of all
operator algebras with morphisms defined by
representations (up to unitary equivalence) of algebras. The universal
repelling object of the category~$\mathcal G$ is the algebra \algebra{G} of
group
convolutions on $\object{Exp}\algebra{g}$.
\end{lem}

Besides a convolution-like calculus, there is another well developed
type of calculus, namely, the calculus of pseudodifferential
operators (PDO), which is highly useful in analysis.
 Let us remind, a
PDO $\object{Op}a(x,\xi)$~\cite{Hormander85,Shubin87,MTaylor81,HWeyl}
with the {\em
Weyl symbol\/} $a(x,\xi)$ is defined by the formula:
\begin{equation}\label{eq:pdo}
[\object{Op}a](x,\xi)u(y)=\int_{\Space{R}{n}\times\Space{R}{n}}
a(\frac{x+y}{2},\xi)\,e^{i\scalar{y-x}{\xi}}u(x)\,dx\,d\xi.
\end{equation}
It was shown that relative convolutions are PDOs for many types studed
before (for example, families of
convolutions~\cite[Proposition~1.1]{MTaylor84} and meta-Heisenberg
group~\cite{Folland94}). But PDOs itself are relative convolutions
induced by the Heisenberg group (see Example~\ref{ex:pdo}). Thus if we
consider morphisms at categories of relative convolution algebras only
up to smooth operators, then we
have\footnote{``The Heisenberg group {\ldots} is basic to this paper
(and much of the rest of the word)''~\cite{Howe80b}.}
\begin{thm}\label{th:pdo} The category $\mathcal G$ of relative convolution
algebras induced by a Lie algebra \algebra{g}, $\dim \algebra{g}= N$ is
(up to smooth operators) a sub-category of the category ${\mathcal H}^N$ of
relative convolution algebras induced by the Heisenberg group
\Heisen{N}.
\end{thm}
\begin{proof} It is well known~\cite[Proposition~1.1]{MTaylor84}, that
group convolutions on $\object{Exp}\algebra{g}$ are (up to smooth
operators, at least) PDOs on \Space{R}{N} and thus the algebra of
convolutions belongs to ${\mathcal H}^N$. The rest
of the assertion is given by Theorem~\ref{th:nature} and
Lemma~\ref{le:category}.
\end{proof}

Of course, there is still a huge distance between Theorem~\ref{th:pdo}
and the real PDO-like calculus of relative convolutions.

Calculus of PDO is a representation~\cite{Howe80b} of
a group convolution calculus on the Heisenberg group (the simplest
nilpotent Lie group). Thus Theorem~\ref{th:pdo} establishes
connection between relative convolutions and convolutions on the
Heisenberg group. Therefore it is not surprising, that for convolutions
on nilpotent Lie groups this connection has a very simple
form~(see~\cite[Theorem~1]{Kisil93e}).

\section{Applications to Complex and Hypercomplex
Analysis}\label{se:analysis}
In this and the next Sections we would like to introduce a series of
essential Examples,
which show possible
applications\footnote{``The most interesting aspect of the {\ldots}
theory has to do with the application of this machinery to concrete
examples''~\cite{Guillemin84}.}  of relative convolutions. Due to the
wide
spectrum of Examples only brief descriptions will be given here. We are
going to
return to these subjects in the future papers in this series after an
appropriate development of general theory of relative convolutions.

As was pointed by H\"ormander in~\cite{Anderson69},
formula~\eqref{eq:relative} is closely connected with a partial
differential equation with operator coefficients. Thus it is not
surprizing, that most Examples within this Section are touching
some spaces of holomorphic functions, which are solutions to
corresponding
equations. In Example~\ref{ex:clifford} this connection will be used
directly.
\begin{example}\label{ex:pdo} Let us consider a combination of
Examples~\ref{ex:euclid-conv} and~\ref{ex:multiplication}. Namely, let
\algebra{g} has a frame consisting from vector fields
\begin{equation}\label{eq:pdo-frame}
X_j=y_j,\ X_j^e=\frac{1}{i}\frac{\partial }{\partial y_j},\ 1\leq
j\leq N,
\end{equation}
which operate on $S=\Space{R}{N}$ by the obvious way. Note, that these
vector
fields have exactly the same commutators~\eqref{eq:heisen-comm} as left
(right) fields from~\eqref{eq:frame-left} (or~\eqref{eq:frame-right}) if
we put $X_0^{l(r)}=iI$. Then an operator of relative convolution with
the kernel $\widehat{k}(x,\xi),\ x,\xi\in\Space{R}{N}$ has form
(see~\cite[\S~1.3]{MTaylor86})
\begin{eqnarray*}
[Kf](y)&=&(2\pi)^{-N}\int_{\Space{R}{2N}} \widehat{k}(x,\xi)\,
           e^{i(\sum_1^Nx_jy_j-\sum_1^N \xi\frac{\partial }{i\partial
y_j})}f(y)
           \,dx \,d\xi\\
&=&(2\pi)^{-N}\int_{\Space{R}{2N}} k(\frac{x+y}{2},\xi)\,
           e^{i(y-x)\xi} f(x)\,dx\,d\xi,
\end{eqnarray*}
i.~e.~it exactly defines the {\em Weyl functional calculus\/} (or the
Weyl
quantization)
\cite{Hormander85,Shubin87,MTaylor81,HWeyl}. This forms a very
important tool for the theory of differential equations and quantum
mechanics.

Our definition of relative convolution operators obviously generalizes
the
Weyl functional calculus from the case of the Heisenberg group to an
arbitrary exponential Lie group. A deep investigation of the role of the
Heisenberg group (and its different representations) at PDO calculus and
harmonic analysis on real line may be found at~\cite{Howe80a,Howe80b}.
Unfortunately, usual  definitions of PDO symbol classes $S^m$ destroy
the natural
symmetry between $x$ and $\xi$ and this restricts applications of
harmonic
analysis to the (standard) PDO theory. Particulary, $S^m(\Space{R}{2n}
)$ is not invariant under all symplectic morphisms of \Space{R}{2n},
which
are induced by automorphisms of the Heisenberg group (see
Corollary~\ref{co:auto}).
\end{example}
\begin{example}\label{ex:toeplitz} Now we follow the
paper~\cite{BoutetGuill85} footsteps but only in the original context of
complex analysis (see also an elegant survey
in~\cite[Chap.~XII]{Stein93}).
Let $\Space{U}{n}$ be an {\em upper half-space\/}
\begin{displaymath}
\Space{U}{n}=\{z\in\Space{C}{n+1}\such \object{Im}z_{n+1} >
\sum_{j+1}^n\modulus{z_j}^2\},
\end{displaymath}
which is a domain of holomorphy of functions of $n+1$ complex variables.
Its boundary
\begin{displaymath}
b\Space{U}{n}=\{z\in\Space{C}{n+1}\such \object{Im}z_{n+1} =
\sum_{j+1}^n\modulus{z_j}^2\}
\end{displaymath}
may be naturally identified with the Heisenberg group \Space{H}{n}.
 One can introduce the {\em Szeg\"o projector\/} $R$ as the orthogonal
projection of $\FSpace{L}{2}(\Heisen{n})$ onto its subspace
$\FSpace{H}{2}(\Heisen{n})$ (the {\em Hardy space\/}) of boundary values
of holomorphic functions on the upper half-space \Space{U}{n}. Then a
{\em Toeplitz operator\/}~\cite{BoutetGuill85}
on \Heisen{n} with the pre-symbol $Q$ is an operator of the form $T_Q=RQR $
where $Q: \FSpace{L}{2}(\Heisen{n}) \rightarrow
\FSpace{L}{2}(\Heisen{n})$ is a pseudodifferential operator. Obviously
$T_Q: \FSpace{H}{2}(\Heisen{n}) \rightarrow \FSpace{H}{2}(\Heisen{n})$.
The invariance of the {\em tangential Cauchy--Riemann equations\/} under
right shifts of
\Heisen{n} implies that the Szeg\"o projector can be realized as a
(left)
convolution operator on \Heisen{n}~\cite{Gindikin64} (see also
Corollary~\ref{co:szego}). Thus the algebra of
the Toeplitz
operators on  \Heisen{n} can be naturally imbedded into the algebra of
(pseudodifferential) operators generated by left group convolutions on
\Heisen{n} and PDO.

First, let us consider a case of a pre-symbol $Q$ taken from usual
Euclidean
convolutions on $\Heisen{n}\cong\Space{R}{2n+1}$. Left convolutions on
\Heisen{n}  are generated by
vector fields $X_j^l$ from~\eqref{eq:frame-left} and Euclidean
convolutions are induced by fields (see Example~\ref{ex:euclid-conv})
\begin{displaymath}
X_j^e=\frac{1}{i}\frac{\partial }{\partial y_j},\ 1\leq j\leq N.
\end{displaymath}
But two frames of vector fields $\{X_j^l,X_j^e\}$ and $\{X_j^l,X_j^r\}$
define
just the same action on \Heisen{n}. Therefore we have an embedding of
the Toeplitz operators with Euclidean convolution pre-symbols to the
algebra of
two-sided convolutions on \Heisen{n} from Example~\ref{ex:two-sided}.

If we now allow $Q$ to be a general PDO from Example~\ref{ex:pdo}, then
we should consider a joint operation of vector fields $X_j,\ X_j^e$
from~\eqref{eq:pdo-frame} and $X_j^l$ from~\eqref{eq:frame-left}. Again
one can pass to equivalent frame defined by $X_j,\ X_j^l,\ X_j^r$.
The algebra of relative convolutions defined by the last frame is
the algebra of
operators generated by two-sided convolutions on \Heisen{n}  and operators of
multiplication by functions, which form a
meta-Heisenberg group~\cite{Folland94}. Such algebras for continuous
multipliers were studied in~\cite{Kisil93e,Kisil94a}.

There is our conclusion\footnote{This is an answer to the reasonable
question of E.~Stein:``Does the algebra of two-sided convolutions
contain at least one interesting operator?''}:
\begin{prop}
The algebra of the Toeplitz operators with PDO
pre-symbols is naturally imbedded into the algebra of relative
convolutions
generated by two-sided convolutions on \Heisen{n} and operators of
multiplication by functions.
\end{prop}
For the first time the algebra of the Toeplitz
operators with two-sided convolution pre-symbols was studied
at~\cite{Kisil91}.
\end{example}
We are going to consider another problem from complex analysis.
\begin{example}\label{ex:fock}
Let $\FSpace{L}{2}(\Space{C}{n},d\mu_{n})$ be a space of all
square-integrable functions on
$\Space{C}{n}$ with respect to the Gaussian measure
\begin{displaymath}
d\mu_{n}(z)=\pi^{-n}e^{-z\cdot\overline{z}}dv(z),
\end{displaymath}
where $dv(z)=dx\,dy$ is the usual Euclidean volume measure on
$\Space{C}{n}=\Space{R}{2n}$. Denote
by $P_{n}$ the orthogonal Bargmann projector of
$\FSpace{L}{2}(\Space{C}{n},d\mu_{n})$ onto the {\em Segal--Bargmann\/} or {\em
Fock space\/}
$\FSpace{F}{2}(\Space{C}{n})$, namely, the subspace of
$\FSpace{L}{2}(\Space{C}{n},d\mu_{n})$ consisting of all entire
functions. The Fock space $\FSpace{F}{2}(\Space{C}{n})$ was introduced
by Fock~\cite{Fock32} to give an alternative representation of the
Heisenberg group in quantum mechanics. The rigorous  theory of
$\FSpace{F}{2}(\Space{C}{n})$ was developed by
Bargmann~\cite{Bargmann61} and Segal~\cite{Segal60}. Such a theory is
closely connected with representations of the Heisenberg group (see
also~\cite{Folland89,Guillemin84,Howe80b}), but studies of the Bargmann
projector and the
associated Toeplitz operators are
usually based on the Hilbert spaces technique\footnote{Or, at least, do not use
harmonic analysis
directly.}. There is a strong reason
for
this: the Bargmann projection $P_n$ is not a {\em group\/} convolution
for any group. However, it is possible to consider $P_n$ as a {\em
relative\/} convolution.

Let us consider the group of Euclidean shifts $a: z\mapsto z+a$ of
\Space{C}{n}. To make unitary operators on
$\FSpace{L}{2}(\Space{C}{n},d\mu)$  from the shifts we should multiply by
the special weight function:
\begin{equation}\label{eq:fock-shift}
a: f(z)\mapsto f(z+a)e^{-z\bar{a}-a\bar{a}/2}.
\end{equation}
 It is obvious, that~\eqref{eq:fock-shift} defines a unitary
representation~\cite{Howe80b} of the $(2n+1)$-dimensional Heisenberg
group on
$\FSpace{L}{2}(\Space{C}{n},d\mu)$, which preserves the Fock space
$\FSpace{F}{2}(\Space{C}{n})$. Thereafter all
operators~\eqref{eq:fock-shift} should commute with $P_n$. Unitary
operators~\eqref{eq:fock-shift} have such infinitesimal generators
\begin{displaymath}
i \sum_{k=1}^{n}\left( a_j'(\frac{\partial }{\partial z_j'}-z_j'-
iz_j'')+
a_j''(\frac{\partial }{\partial z_j''}-z_j''+iz_j')\right),
\end{displaymath}
where $a=(a_1,\ldots,a_n), z=(z_1,\ldots,z_n)\in\Space{C}{n}$ and
$a_j=(a_j',a_j''), z=(z_j',z_j'')\in\Space{R}{2}$. This linear space of
generators has the frame
\begin{equation}\label{eq:fshift-frame}
A^{f\prime}_j=\frac{1}{i}\left(\frac{\partial }{\partial z_j'}-z_j'-
iz_j''\right),\
A^{f\prime\prime}= \frac{1}{i}\left(\frac{\partial }{\partial z_j''}-
z_j''+iz_j'\right).
\end{equation}
Operators~\eqref{eq:fshift-frame} still should commute with the Bargmann
projector $P_n$
and we can expect that $P_n$ should be a relative convolution with
respect to operators
\begin{equation}\label{eq:fconv-frame}
X^{f\prime}_j=\frac{1}{i}\left(\frac{\partial }{\partial z_j'}-
z_j'+iz_j''\right),\
X^{f\prime\prime}= \frac{1}{i}\left(\frac{\partial }{\partial z_j''}-
z_j''-iz_j'\right),
\end{equation}
which commute with all operators~\eqref{eq:fshift-frame} and together
with operator $2I$ form a self-adjoint representation of
$\algebra{h}_n$. Indeed, we have
\begin{prop}\label{pr:bargmann} The Bargmann projector is a relative
convolution induced by
the Weyl-Heisenberg Lie algebra $\algebra{h}_n$, which have an operation on
\Space{C}{n}
defined by~\eqref{eq:fconv-frame}. Its kernel $b(t,\zeta),\ t\in
\Space{R}{},\ \zeta\in\Space{C}{n}$ is defined by the formula:
\begin{displaymath}
\widehat{b}(t,\zeta)=2^{n+1/2} e^{-1} e^{-
(t^2+\zeta\bar{\zeta}/2)}.
\end{displaymath}
\end{prop}
\begin{proof} One can make an easy exercise with integral transforms:
\begin{eqnarray}
[P_nf](z)&=&(2\pi)^{-n-1/2}\int_{\Space{R}{}} \int_{\Space{C}{n}} 2^{n+1/2}\,
e^{-(t^2+1+\zeta\bar{\zeta})/2} \nonumber \\
     &&\eqnright\times e^{-i(t\cdot 2I+\sum_{k=1}^n
(\zeta_j'X_j^{f\prime}+\zeta_j''X_j^{f\prime\prime}))}\,
f(z)\,dt\,d\zeta\nonumber \\
    &=&\pi^{-n-1/2}\int_{\Space{R}{}}e^{-(t^2+1+2it)}
\,dt\,\int_{\Space{C}{n}}\,  e^{-(\zeta\bar{\zeta})/2}\, e^{-i\sum_{k=1}^n
(\zeta_j'X_j^{f\prime}+\zeta_j''X_j^{f\prime\prime})}\,
f(z)\,d\zeta\nonumber \\
    &=&\pi^{-n}\int_{\Space{C}{n}} e^{-\zeta\bar{\zeta}/2}\, e^{-
i\sum_{k=1}^n
(\zeta_j'X_j^{f\prime}+\zeta_j''X_j^{f\prime\prime})}\,
f(z)\,d\zeta\nonumber \\
     &=&\pi^{-n}\int_{\Space{C}{n}} e^{-\zeta\bar{\zeta}/2}\, e^{-
\sum_{k=1}^n
(\zeta_j'\frac{\partial }{\partial z_j'}
+\zeta_j''\frac{\partial }{\partial z_j''}- \zeta_j \bar{z}_j)}\,
f(z)\,d\zeta\nonumber \\
   &=&\pi^{-n}\int_{\Space{C}{n}} e^{-\zeta\bar{\zeta}/2}\, e^{-
\zeta\bar{\zeta}/2}\,
e^{\zeta \bar{z}}\, e^{-\sum_{k=1}^n(\zeta_j'\frac{\partial }{\partial
z_j'} +\zeta_j''\frac{\partial }{\partial z_j''})}\, f(z)\,d\zeta\nonumber
\\
   &=&\pi^{-n}\int_{\Space{C}{n}} e^{(\bar{z}-\bar{\zeta})\zeta}\, f(z-
\zeta)\,
d\zeta\nonumber \\
   &=&\pi^{-n}\int_{\Space{C}{n}}\, e^{\bar{w}(z-w)}\,
f(w)\,dw\label{eq:bargmann},
\end{eqnarray}
where $\zeta=(\zeta',\zeta'')\in\Space{R}{n} \times \Space{R}{n},\ w=z-
\zeta$. Formula~\eqref{eq:bargmann} is the well known expression
for the Bargmann projector~\cite{Bargmann61}.
\end{proof}
See also Corollary~\ref{co:bargmann} for an alternative proof.

Following to
papers~\cite{BergCob86,BergCob87,Coburn90} we now consider
the Toeplitz operator of the form $T_a=P_n a(z)I$, where $a(z)I$ is an
operator of multiplication by a function $a(z)$. For the previous
reasons we can handle them as relative convolutions generated by
operators $X^f_j$ from~\eqref{eq:fconv-frame} and operators $X_j=z_j I$.
We can easily describe all non-zero commutators:
\begin{equation}\label{eq:fock-comm}
[X^{f\prime}_j,X^{f\prime\prime}_j]=2iI,\ [X^{f\prime}_j,X'_j]=iI,\
[X^{f\prime\prime}_j,X''_j]=iI.
\end{equation}
So they form  a $(4n+1)$-dimensional nilpotent step 2 Lie algebra.
Particulary sub-algebra spanned on the vectors $X^{f\prime}_j,\
X^{f\prime\prime}_j$ and $iI$ is isomorphic\footnote{This explains why
``a
critical ingredient in our analysis is an averaging operation over the
Segal--Bargmann representation of the Heisenberg
group''~\cite{BergCob87}.} to our
constant companion $\algebra{h}_n$.
\begin{prop} There is a natural embedding of the Toeplitz operator
algebra on
the Fock space into the algebra of relative convolutions induced by the
Lie
algebra with commutation relations~\eqref{eq:fock-comm}.
\end{prop}
We will continue our discussion of Toeplitz operators on an ``abstract
nonsense''\cite{Howe80b} level at Section~\ref{se:coherent}.
\end{example}
\begin{example}\label{ex:clifford}
Let now $X_j$ be generators of the Clifford algebra \Cliff{0}{n} (we use
book~\cite{DelSomSou92} as a standard reference within this Example).
This
means that the following {\em anti-commutation\/} relations hold
(compare
with~\eqref{eq:heisen-comm}):
\begin{equation}\label{eq:anti-comm}
\{X_i,X_j\}:=X_i X_j+X_j X_i=-2\delta_{ij}X_0,
\end{equation}
where $X_0=I$. Function $f:\Space{R}{n}\rightarrow\Cliff{0}{n}$ is
called {\em monogenic\/} if it satisfies the {\em Dirac equation\/}
\begin{equation}\label{eq:dirac}
Df:=\frac{\partial f(y)}{\partial y_0} -\sum_{j=1}^n X_j\frac{\partial
f(y)}{\partial y_j}=0\ \mbox{ or }\ \frac{\partial f(y)}{\partial y_0}
=\sum_{j=1}^n X_j\frac{\partial
f(y)}{\partial y_j}.
\end{equation}
The success of Clifford analysis is mainly explained because
the Dirac operator~\eqref{eq:dirac} factorizes the Laplace operator
$\Delta=\sum_0^n\frac{\partial ^2}{\partial x_j^2}$.

H\"ormander's remark from paper~\cite{Anderson69} gives us by
the
fundamental solution to the Dirac equation in the form
\begin{eqnarray*}
K(y)&=&\fourier{\eta\rightarrow y} e^{-iy_0\sum_{j=1}^n\eta_jX_j}\\
    &=&\int_{\Space{R}{n}}e^{i\sum_{j=1}^n y_j\eta_j}\, e^{-
iy_0\sum_{j=1}^n\eta_jX_j}\, d\eta.
\end{eqnarray*}
(``Simply take the Fourier transform with respect to the spatial
variables, and solve the equation in
$y_0$''~\cite{Anderson69}). Otherwise, any solution $f(y)$
to~\eqref{eq:dirac} is given by a
convolution of some function $\widetilde{f}(y)$ on \Space{R}{n-1}  and
the fundamental
solution $K(y)$. On the contrary, a convolution $K(y)$ with any function
is a solution to~\eqref{eq:dirac}. We have:
\begin{eqnarray}
[K*f](y)&=&\int_{\Space{R}{n}} K(y-t)\,f(t)\,dt\nonumber \\
        &=&\int_{\Space{R}{n}}\int_{\Space{R}{n}} e^{-i\sum_{j=1}^n (y_j-
t_j)\eta_j}\, e^{-iy_0\sum_{j=1}^n\eta_jX_j}\, d\eta f(t)\,dt\nonumber \\
        &=&\int_{\Space{R}{n}}\int_{\Space{R}{n}}e^{i\sum_{j=1}^n t_j\eta_j}\,
e^{-i\sum_{j=1}^n\eta_j(y_0 X_j-y_j X_0)}\, d\eta f(t)\,dt\nonumber \\
        &=&\int_{\Space{R}{n}} e^{-i\sum_{j=1}^n\eta_j(y_0 X_j-y_j
X_0)}\int_{\Space{R}{n}} e^{i\sum_{j=1}^n t_j\eta_j}\, f(t)\,dt\, d\eta
\nonumber \\
        &=&\int_{\Space{R}{n}} e^{-i\sum_{j=1}^n\eta_j(y_0 X_j-y_j X_0)}\,
\widehat{f}(-\eta)\, d\eta. \label{eq:malonek}
\end{eqnarray}
Equation~\eqref{eq:malonek} defines relative
convolution~\eqref{eq:relative} with the kernel $f(-y)$ and the Lie
algebra
of vector fields
\begin{equation}\label{eq:monom}
\{y_0 X_j-y_j X_0\},\ 1\leq j \leq n.
\end{equation}
Thus, at least formally, any solution to the Dirac
equation~\eqref{eq:dirac} can be written as a function (see
Remark~\ref{re:function}) of $n-1$
monomials~\eqref{eq:monom}:
\begin{equation}\label{eq:laville}
\breve{f}(y_0,y_1,\ldots,y_n)= f(y_0 X_1-y_1 X_0, y_0 X_2-y_2 X_0,\ldots,y_0
X_n-y_n X_0).
\end{equation}
Another significant remark: if we fix the value $y_0=0$ in~\eqref{eq:laville}
we easily obtain: $\breve{f}(0,y_1,\ldots,y_n)= f(-y_1 X_0, -y_2
X_0,\ldots,-y_n X_0)=f(-y_1 , -y_2 ,\ldots,-y_n )$. Thus we may consider the
function $\breve{f}(y_0,y_1,\ldots
y_n)$ of $n+1$ variables as {\em analytical expansion\/} for the function
$f(y_1,\ldots,y_n)$ of $n$ variables (compare with~\cite{Laville91}).

Using the power series decomposition for the exponent one can see that
formula~\eqref{eq:malonek} defines the permutational (symmetric) product
of monomials~\eqref{eq:monom}. The significant role of such monomials
and functions generated by them in Clifford analysis was recently
discovered by Laville~\cite{Laville87} and  Malonek~\cite{Malonek93}. But
during our consideration we
used only the commutation relation $[X_0,X_j]=0$ and never used the
anti-commutation relations~\eqref{eq:anti-comm}. Thus
formula~\eqref{eq:malonek} {\em is true and may be useful without
Clifford analysis\/} as well.
\begin{prop}
Any solution to equation~\eqref{eq:dirac}, where $X_j$ are arbitrary
self-adjoint operators, is given as arbitrary function of $n$
monomials~\eqref{eq:monom} by the formula~\eqref{eq:malonek}.
\end{prop}
It is possible also to introduce the notion of
the {\em differentiability\/}~\cite{Malonek93}  for
solutions to~\eqref{eq:dirac}, namely, an increment of any solution
to~\eqref{eq:dirac} may be approximated up to infinitesimals of the
second
order by a linear function of monomials~\eqref{eq:monom}.

Due to physical application we will consider equation
\begin{equation}\label{eq:mass}
\frac{\partial f}{\partial y_0}=(\sum_{j=1}^n X_j\frac{\partial
}{\partial y_j}+M)f,
\end{equation}
where $X_j$ are arbitrary self-adjoint operators and $M$ is a bounded
operator commuting with all $X_j$.
\begin{rem}
Then $X_j$ are generators~\eqref{eq:anti-comm} of the Clifford algebra
and $M=M_\alpha $ is an operator of
multiplication from the {\em right\/}-hand side by the Clifford number $\alpha
$,
differential operator~\eqref{eq:mass} factorizes the Helmholtz operator
$\Delta +M_{\alpha^2}$. Equation~\eqref{eq:mass} is known in quantum
mechanics as the
{\em Dirac
equation for a particle with a non-zero rest
mass\/}~\cite[\S 20]{BerLif82}, \cite[\S 6.3]{BogShir80} and
\cite{Kravchenko95a}.
\end{rem}
Simple modification of the previous calculations gives us the
following
result
\begin{prop}
Any solution to equation~\eqref{eq:mass}, where $X_j$ are arbitrary
self-adjoint operators and $M$ commutes with them, is given by the
formula
\begin{displaymath}
e^{y_0 M}\int_{\Space{R}{n}} e^{-i\sum_{j=1}^n\eta_j(y_0 X_j-y_j X_0)}
\widehat{f}(-\eta)\, d\eta,
\end{displaymath}
where $f$ is an arbitrary function on \Space{R}{n-1}.
\end{prop}
\begin{rem}\label{re:riesz}
Connection between relative convolutions and Clifford analysis is
two-sided. Not only relative convolutions are helpful in Clifford
analysis but also Clifford analysis may be used for developing
the relative
convolution technique. Indeed, we have defined a relative convolution as
a function of operators $X_j$ representing a Lie algebra \algebra{g}.
To do this we have used the Weyl function calculus
from~\cite{Anderson69}. Meanwhile, for a pair of self-adjoint operators
$X_1,\ X_2$ the alternative Riesz calculus~\cite[Chap.~XI]{RieszNagy55}
is given by the formula
\begin{equation}\label{eq:riesz}
f(X_1+iX_2)=\int_\Gamma f(\tau)(\tau I-(X_1+iX_2))^{-1}\,d\tau.
\end{equation}
As shown in~\cite{Anderson69}, these two calculuses are essentially the
same
in the case of a pair of bounded operators. To extend the Riesz calculus
for arbitrary $n$-tuples of bounded operators $\{X_j\}$ it seems natural
to use Clifford analysis (see, for example,~\cite{DelSomSou92}), which
is an analogy to one-dimensional complex analysis. Then one can define a
function of arbitrary $n$-tuple of bounded self-adjoint operators $\{X_j\}$ in
such a way (compare with \eqref{eq:riesz}):
\begin{displaymath}
f(X_1,X_2,\ldots,X_n)=\int_\Gamma
f(\tau_1,\tau_2,\ldots,\tau_n)\,K(\tau_1I-X_1,\tau_2I-X_2,\ldots,\tau_nI-
X_n)\,d\tau,
\end{displaymath}
where $K(\tau_1,\tau_2,\ldots,\tau_n)$ is the Cauchy kernel fron the Clifford
analysis~\cite{DelSomSou92} .
\end{rem}
\end{example}

\section{Coherent States}\label{se:coherent}
This Section is a bridge between the previous and the next ones. Here we
will give a new glance on some constructions of
Examples~\ref{ex:toeplitz} and~\ref{ex:fock}. We also provide a
foundation for the further investigation in Section~\ref{se:physics}.

\subsection{General Consideration}\label{ss:general}
 Coherent states are a useful tool in quantum theory and
have a lot of essentially different definitions~\cite{Perelomov86}.
Particulary, they were described by Berezin
in~\cite{Berezin72,Berezin86,Guillemin84} concerning so-called
covariant and
contravariant (or Wick and anti-Wick) symbols of operators
(quantization).
\begin{defn} They say that the Hilbert space $H$ has a system of
coherent states $\{f_\alpha \}, \alpha \in G$ if for any $f\in H$
\begin{equation}\label{eq:norm}
\scalar{f}{f}=\int_G \modulus{\scalar{f}{f_\alpha }}^2\, d\mu.
\end{equation}
\end{defn}
 This definition does not take in account that within
coherent states  a group structure frequently occurs and is always
useful~\cite{Perelomov86}.
For example, the original consideration of Berezin is connected with the
Fock space, where coherent states are functions $e^{za-a\bar{a}/2}$. The
representation of the Heisenberg group on the Fock space was exploited
in Example~\ref{ex:fock}. Another type of coherent states with a group
structure is given by the vacuum vector and operators of creation and
annihilation, which represents group \Space{Z}{}. Thus we would like to
give an alternative definition.

\begin{defn}\label{de:coherent1} We will say that the Hilbert space $H$
has a system of
coherent states $\{f_g \},\ g \in G$ if
\begin{enumerate}
\item\label{it:group} There is a representation $T: g\mapsto T_g$ of the
group
$G$ by unitary operators $T_g$ on $H$.
\item\label{it:vector} There is a vector $f_0\in H$ that for $f_g=T_g
f_0$ and arbitrary $f\in H$ we have
\begin{equation}\label{eq:norm1}
\scalar{f}{f}=\int_G \modulus{\scalar{f}{f_g }}^2\, d\mu,
\end{equation}
where integration is taken over the Haar measure $d\mu$ on $G$.
\end{enumerate}
\end{defn}
Because this construction independently arose in different contexts,
vector $f_0$ has many various names: the {\em vacuum vector\/}, the {\em ground
state\/}, the
{\em mother wavelet\/} etc. Modifications of
definition~\ref{de:coherent1} for other cases are
discussed~\cite[\S~2.1]{Perelomov86}.
Equation~\eqref{eq:norm1} implies, that vector $f_0$ is a cyclic vector
of the representation $G$ on $H$.
\begin{lem}
Let $T$ is a unitary representation of a group $G$ in a Hilbert space
$H$. Then there exists such $f_0\in H$ that equality~\eqref{eq:norm1} holds.
Moreover, if the representation $T$ is irreducible, then one can take an
arbitrary non-zero vector of $H$ (up to a scalar factor) as the vacuum vector.
\end{lem}
\begin{proof}
Let us fix some Haar measure $d\mu$ on $G$. (Different Haar measures
are different on a scalar factor). If the representation $T$ is irreducible,
then an arbitrary vector $f\in H$ is cyclic and we may put $f_0=c^{-
1/2}f$, where
\begin{displaymath}
c=\frac{\int_G \modulus{\scalar{f}{T_g f}}^2\, d\mu(g)}{\scalar{f}{f}}.
\end{displaymath}
It is easy to verify, that for the $f_0$ equality~\eqref{eq:norm1}
holds~\cite[\S~2.3]{Perelomov86}.

Let now $T$ is an arbitrary representation, then we can
decompose~\cite[\S~8.4]{Kirillov76} it onto a direct
integral~\cite[\S~10]{Dixmier69} of irreducible representations
$T=\int_Y T_\alpha\, d\alpha $ on the space $H=\int_Y H_\alpha\,
d\alpha$. Again we can take an arbitrary vector $f=\int_Y
f_\alpha\, d\alpha \in H=\int_Y H_\alpha\, d\alpha$, such that $f_\alpha $ are
non-zero almost everywhere, and put
$f_0=\int_Y f_\alpha c^{-1/2}_\alpha d\alpha$, where
\begin{displaymath}
c_\alpha =\frac{\int_G \modulus{\scalar{f_\alpha }{T_g f_\alpha }_\alpha
}^2 d\mu(g)}{\scalar{f_\alpha }{f_\alpha }_\alpha }.
\end{displaymath}
In view of $\norm{f}=\int_Y \norm{f_\alpha }\,d\alpha $, the proof is
complete.
\end{proof}
By the way, a polarization
of~\eqref{eq:norm1} gives us the equality
\begin{equation}\label{eq:product}
\scalar{f_1}{f_2}=\int_G \scalar{f_1}{f_g }
\overline{\scalar{f_2}{f_g
}}\, d\mu.
\end{equation}
Thus we have an isometrical embedding $E: H\rightarrow
\FSpace{L}{2}(G,d\mu)$ defined by the formula
\begin{equation}\label{eq:embedding}
E: f \mapsto f(g)=\scalar{f}{f_g}=\scalar{f}{T_g f_0}= \scalar{T^*_g
f}{f_0}= \scalar{T_{g^{-1}}f}{f_0}.
\end{equation}
We will consider $\FSpace{L}{2}(G,d\mu)$ both as a linear space of
functions and as an operator algebra with respect to the left and
right group
convolution operations:
\begin{eqnarray}
{}[f_1*f_2]_l(h)&=&\int_G f_1(g)f_2(g^{-1}h)\, d\mu(g),
\label{eq:conv-l}\\
{}[f_1*f_2]_r(h)&=&\int_G f_1(g)f_2(hg)\, d\mu(g).\label{eq:conv-r}
\end{eqnarray}
For a simplicity we will assume, that $G$ is unimodular (the left
and the right Haar measures on $G$ coincide) and that
$\FSpace{L}{2}(G,d\mu)$ is closed under the group convolution.
Thus the construction under consideration may be regarded as a natural
embedding
of the linear space $H$ to the operator algebra $\bos(H)$.

Let $\FSpace{H}{2}(G,d\mu)\subset\FSpace{L}{2}(G,d\mu)$ will denote the
image of $H$ under embedding $E$. It is clear, that
$\FSpace{H}{2}(G,d\mu)$ is a liner subspace of $\FSpace{L}{2}(G,d\mu)$,
which does not coincide with the whole $\FSpace{L}{2}(G,d\mu)$ in general.
One can see, that
\begin{lem}\label{le:h2-invar}
Space $\FSpace{H}{2}(G,d\mu)$ is invariant under left shifts
on $G$.
\end{lem}
\begin{proof} Indeed, for every $f(g)\in \FSpace{H}{2}(G,d\mu)$
the function
\begin{displaymath}
f(h^{-1}g)=\scalar{f}{T_{h^{-1}g}f_0}=
\scalar{f}{T_{h^{-1}}T_gf_0}=\scalar{T_{h}f}{T_gf_0}=[T_{h}f](g)
\end{displaymath}
also belongs to $\FSpace{H}{2}(G,d\mu)$.
\end{proof}
If
$P:\FSpace{L}{2}(G,d\mu)\rightarrow \FSpace{H}{2}(G,d\mu)$ the
orthogonal projector on $\FSpace{H}{2}(G,d\mu)$, then due to
Lemma~\ref{le:h2-invar} it should commute
with all left shifts and thus we get immediately
\begin{cor}
Projector $P:\FSpace{L}{2}(G,d\mu)\rightarrow \FSpace{H}{2}(G,d\mu)$ is
a right convolution on $G$.
\end{cor}
The following Lemma characterizes {\em linear
subspaces\/} of $\FSpace{L}{2}(G,d\mu)$ invariant under shifts in the term of
{\em convolution
algebra\/} $\FSpace{L}{2}(G,d\mu)$ and seems to be of the separate
interest.
\begin{lem} \label{le:ideal}
A closed linear subspace $H$ of $\FSpace{L}{2}(G,d\mu)$ is invariant
under left (right) shifts if and only if $H$ is a left (right) ideal of
the right group convolution algebra $\FSpace{L}{2}(G,d\mu)$.

A closed linear subspace $H$ of $\FSpace{L}{2}(G,d\mu)$ is invariant
under left (right) shifts if and only if $H$ is a right (left) ideal of
the left group convolution algebra $\FSpace{L}{2}(G,d\mu)$.
\end{lem}
\begin{proof} Of course we consider only the ``right-invariance and
right-convolution'' case. Then the other three
cases are analogous. Let $H$ be a closed linear subspace of
$\FSpace{L}{2}(G,d\mu)$ invariant under right shifts and $k(g)\in H$. We
will show the inclusion
\begin{equation}\label{eq:rt-convol}
[f*k]_r(h)=\int_G f(g)k(hg)\,d\mu(g)\in H,
\end{equation}
for any $f\in\FSpace{L}{2}(G,d\mu)$. Indeed, we can treat
integral~\eqref{eq:rt-convol} as a limit of sums
\begin{equation}\label{eq:rt-sum}
\sum_{j=1}^{N} f(g_j)k(hg_j)\Delta_j.
\end{equation}
But the last sum is simply a linear combination of vectors $k(hg_j)\in
H$ (by the invariance of $H$) with coefficients $f(g_j)$. Therefore
sum~\eqref{eq:rt-sum} belongs to $H$ and this is true for
integral~\eqref{eq:rt-convol} by the closeness of $H$.

Otherwise, let $H$ be a right ideal in the group convolution algebra
$\FSpace{L}{2}(G,d\mu)$ and let $\phi_j(g)\in\FSpace{L}{2}(G,d\mu)$ be
an approximate unit of the algebra~\cite[\S~13.2]{Dixmier69}, i.~e. for
any $f\in\FSpace{L}{2}(G,d\mu)$ we have
\begin{displaymath}
[\phi_j*f]_r(h)=\int_G \phi_j(g)f(hg)\, d\mu(g) \rightarrow f(h)\mbox{,
when } j\rightarrow\infty.
\end{displaymath}
Then for $k(g)\in H$ and for any $h'\in G$ the right convolution
\begin{displaymath}
[\phi_j*k]_r(hh')=\int_G \phi_j(g)k(hh'g)\, d\mu(g)= \int_G
\phi_j(h'^{-1}g')k(hg')\, d\mu(g'),\ g'=h'g,
\end{displaymath}
from the first expression is tensing to $k(hh')$ and from the second
one belongs to $H$ (as a right ideal). Again the closeness of $H$
implies $k(hh')\in H$ that proves the assertion.
\end{proof}
\begin{lem}\label{le:identity}\upshape{(The reproducing property)} For
any $f(g)\in\FSpace{H}{2}(G,d\mu)$ we have
\begin{eqnarray}
{}[f*f_0]_l(g)&=&f(g) \label{eq:l-reproduce}\\
{} [\bar{f}_0*f]_r(g)&=&f(g),\label{eq:r-reproduce}
\end{eqnarray}
where $f_0(g)=\scalar{f_0}{T_gf_0}$ is the function
corresponding to the vacuum vector $f_0\in H$.
\end{lem}
\begin{proof} We again check only the left case and this is just a
simple calculation:
\begin{eqnarray*}
[f*f_0]_l(h)&=& \int_G f(g)\,f_0(g^{-1}h)\,d\mu(g)\\
          &=& \int_G f(g)\scalar{f_0}{T_{g^{-1}h}f_0}\,d\mu(g)\\
          &=& \int_G \scalar{f}{T_g f_0} \,\scalar{f_0}{T_{g^{-1}} T_h
f_0}\,
d\mu(g)\\
          &=& \int_G \scalar{f}{T_g f_0}\, \scalar{T_g f_0}{ T_h f_0}\,
d\mu(g)\\
          &=& \int_G \scalar{f}{T_g f_0}\, \overline{\scalar{ T_h f_0}{T_g
f_0}}\,
d\mu(g)\\
          &\stackrel{(*)}{=}& \scalar{f}{T_h f_0}\\
    &=&f(h).
\end{eqnarray*}
Here transformation $(*)$ is based on~\eqref{eq:product} and we have
used the unitary property of the representation $T$.
\end{proof}
The following general Theorem easily follows from the previous Lemmas
\begin{thm}\label{th:projector}
The orthogonal projector $P: \FSpace{L}{2}(G,d\mu) \rightarrow
\FSpace{H}{2}(G,d\mu)$ is a right convolution on $G$ with the kernel
$\bar{f}_0(g)$ defined by the vacuum vector.

\textup{To put in the Archimedes-like words,} let me a representation of
group $G$ on $H$ with a cyclic vector $f_0$ and I
will construct the projector $P: \FSpace{L}{2}(G,d\mu) \rightarrow
\FSpace{H}{2}(G,d\mu)\cong H$.
\end{thm}
\begin{proof}
Let $P$ be the operator of right convolution~\eqref{eq:conv-r} with the
kernel
$\bar{f}_0(g)$. By the
Lemma~\ref{le:h2-invar} $\FSpace{H}{2}(G,d\mu)$ is an invariant linear
subspace of $\FSpace{L}{2}(G,d\mu)$. Thus by Lemma~\ref{le:ideal} it is
an ideal under convolution operators. Therefore the convolution operator
$P$ with the kernel $\bar{f}_0(g)$ from $\FSpace{H}{2}(G,d\mu)$ has an image
belonging to $\FSpace{H}{2}(G,d\mu)$. But by Lemma~\ref{le:identity}
$P=I$
on $\FSpace{H}{2}(G,d\mu)$, so $P^2=P$ on $\FSpace{L}{2}(G,d\mu)$, i.~e.
$P$ is a projector on $\FSpace{H}{2}(G,d\mu)$.

It is
easy to see, that $f_0(g)$ has the property $f_0(g)=\bar{f}_0(-g)$, thus
$P^*=P$, i.~e. $P$ is orthogonal. It may be shown also in such a manner.
Let $f(g)\in \FSpace{L}{2}(G,d\mu)$ be orthogonal to all functions from
$\FSpace{H}{2}(G,d\mu)$. Particulary $f(g)$ should be orthogonal to
$f_0(h^{-1}g)$ (due to the invariance of $\FSpace{H}{2}(G,d\mu)$) for
any $h\in G$. Then
$P(f)=[f*f_0]_l=0$ and we have shown the orthogonality again. This
completes the proof.
\end{proof}
\begin{rem}
The stated left invariance of $\FSpace{H}{2}(G,d\mu)$ and the
representation of $P$ as a
right group convolution have a useful tie with differential equations.
Really, let $X_j,\ j\leq m$ be left-invariant vector fields (i.~e.
left-invariant differential operators) on $G$. If $X_jf_0\equiv 0$ then
$X_jf= 0$ for any $f\in H$. Thus space $\FSpace{H}{2}(G,d\mu)$ may be
characterized as the
space of solutions to the system of equations $X_jf= 0,\ 1\leq j\leq m$.

Another connected formulation: we can think of $P$ as of an integral
operator with the kernel $K(h,g)=f_0(g^{-1}h)=\overline{f_0(h^{-
1}g)}=\scalar{T_g f_0}{T_h f_0}$, then kernel $K(h,g)$  is an {\em
analytic\/} function of $h$ and {\em anti-analytic\/} of $g$.
\end{rem}
\begin{example}\label{ex:transform}
An important class of applications may be treated as follows. Let we
have a space of functions defined on a domain $\Omega\in\Space{R}{n} $.
Let
we have also a transitive Lie group $G$ of automorphisms of $\Omega$.
Then we can construct a unitary representation $T$ of $G$ on
$\FSpace{L}{2}(\Omega)$ by the formula:
\begin{equation}\label{eq:unit-trans}
T_g: f(x) \mapsto f(g(x))\, J^{1/2}_g(x),\ f(x)\in\FSpace{L}{2}(\Omega),\
g\in G,
\end{equation}
where $J_g(x)$ is the Jacobian of the transformation defined by $g$ at
the point $x$.

If we fix some point $x_0\in \Omega$ then we can identify the
homogeneous space $G/G_{x_0}$ with $\Omega$ (see notation at the
beginning of Subsection~\ref{ss:notation}). Then left-invariant vector
fields on $G$ may be considered as differential operators on $\Omega$
and convolution operators on $G$ as integral operators on $\Omega$. This
is a way giving {\em integral representations for analytic functions\/}.
\end{example}

\subsection{Classical Results}
We would like to show, how abstract Theorem~\ref{th:projector} and
Example~\ref{ex:transform} are connected
with classical results on the Bargmann, Bergman and Szeg\"o projectors
at the Segal-Bargmann (Fock), Bergman and Hardy spaces respectively. We
will start from a trivial example.
\begin{cor} Let $\{\phi_j\},\ -\infty<j<\infty$ be an orthonormalized
basis of
a Hilbert space $H$. Then
\begin{equation}
B=\sum_{j=-\infty}^\infty \ket{\phi_j}\bra{\phi_j}
\end{equation}
is a reproducing operator, namely, $Bf=f$ for any $f\in H$.
\end{cor}
\begin{proof} We will construct a unitary representation of group
\Space{Z}{} on $H$ by its action on the basis:
\begin{displaymath}
T_k \phi_j = \phi_{j+k},\ k\in \Space{Z}{}.
\end{displaymath}
If we equip \Space{Z}{} with the invariant discreet measure $d\mu(k)=1$
and  select the vacuum vector $f_0=\phi_0$, then all coherent states are
exactly the basis vectors: $f_k=T_k\phi_0=\phi_k$.
Equation~\eqref{eq:product} turns to be exactly the Plancherel formula
\begin{displaymath}
\scalar{f_1}{f_2}=\sum_{j=-\infty}^\infty \scalar{f_1}{T_jf_0}
\overline{\scalar{f_2}{T_jf_0}}=\sum_{j=-\infty}^\infty
\scalar{f_1}{\phi_j}
\overline{\scalar{f_2}{\phi_j}}
\end{displaymath}
and we have obtained the usual isomorphism of Hilbert spaces
$H\cong\FSpace{l}{2}(\Space{Z}{})$ by the formula $f(k)=
\scalar{f}{\phi_{k}}$. Our construction gives
\begin{eqnarray*}
f(k) &=& \sum_{j=-\infty}^\infty \overline{\scalar{f_0}{T_j
f_0}}\scalar{f}{T_{j+k}f_0}\\
     &=& \sum_{j=-\infty}^\infty \delta_{0j} \scalar{f}{\phi_{j+k}}\\
     &=& \scalar{f}{\phi_{k}}.
\end{eqnarray*}
Thus operator $B$ is really identical on $H$. Note, that similar
construction may be given for a case of not orthonormalized frame.
\end{proof}
In spite of simplicity of this construction, it was the (almost) unique
tool to establish of various projectors
(see~\cite[3.1.4]{Rudin80}). Following less trivial Corollaries bring us
back to
Section~\ref{se:analysis}.
\begin{cor}\label{co:bargmann}\textup{\cite{Bargmann61}} The Bargmann projector
on the
Segal--Bargmann space has the kernel
\begin{equation}\label{eq:kr-bargmann}
K(z,w)=e^{\bar{w}(z-w)}.
\end{equation}
\end{cor}
\begin{proof}
Let us define a unitary representation of the Heisenberg group
\Heisen{n} on $\Space{R}{n}$ by the formula~\cite[\S~1.1]{MTaylor86}:
\begin{displaymath}
g=(t,q,p): f(x) \rightarrow T_{(t,q,p)}f(x)=e^{i(2t-\sqrt{2}qx+qp)} f(x-
\sqrt{2}p).
\end{displaymath}
As ``vacuum vector'' we will select the original {\em vacuum vector\/}
$f_0(x)=e^{-x^2/2}$. Then embedding $\FSpace{L}{2}(\Space{R}{n} )
\rightarrow  \FSpace{L}{2}(\Heisen{n} ) $ is given by the formula
\begin{eqnarray}
\widetilde{f}(g)&=&\scalar{f}{T_gf_0}\nonumber \\
     &=&\pi^{-n/4}\int_{\Space{R}{n}} f(x)\,e^{-i(2t-\sqrt{2}qx+qp)}\,e^{-(x-
\sqrt{2}p)^2/2}\,dx\nonumber \\
     &=&e^{-2it-(p^2+q^2)/2}\pi^{-n/4}\int_{\Space{R}{n}} f(x)\,e^{-
((p+iq)^2+x^2)/2+\sqrt{2}(p+iq)x}\,dx\nonumber \\
     &=&e^{-2it-z\bar{z}/2}\pi^{-n/4}\int_{\Space{R}{n}} f(x)\,e^{-
(z^2+x^2)/2+\sqrt{2}zx}\,dx, \label{eq:tr-bargmann}
\end{eqnarray}
where $z=p+iq,\ g=(t,p,q)=(t,z)$. Then $\widetilde{f}(g)$ belong to
$\FSpace{L}{2}(\Heisen{n}, dg)$. It is easy to see, that for every fixed
$t_0$ function $\breve{f}(z)=e^{z\bar{z}/2}\widetilde{f}(t_0,z)$ belongs
to the Segal-Bargmann space, say, is analytic by $z$  and
square-integrable with respect the Gaussian measure $\pi^{-n}e^{-
z\bar{z}}$. The integral in~\eqref{eq:tr-bargmann} is the well known
Bargmann transform~\cite{Bargmann61}. Then the projector
$\FSpace{L}{2}(\Heisen{n},
dg)\rightarrow \FSpace{L}{2}(\Space{R}{n})$ is a convolution on the
\Heisen{n} with the kernel
\begin{eqnarray}
P(t,q,p)&=&\scalar{f_0}{T_gf_0}\nonumber \\
        &=&\pi^{-n/4}\int_{\Space{R}{n}} e^{-x^2/2}\,e^{-i(2t-
\sqrt{2}qx+qp)}\,e^{-(x-\sqrt{2}p)^2/2}\,dx\nonumber \\
        &=&\pi^{-n/4}\int_{\Space{R}{n}} e^{-x^2/2 -i2t+\sqrt{2}iqx-iqp-
x^2/2+\sqrt{2}px -p^2}\,dx\nonumber \\
        &=&e^{-i2t-(p^2+q^2)/2}\pi^{-n/4}\int_{\Space{R}{n}} e^{-(x-
(p+iq)/\sqrt{2})^2}\,dx\nonumber \\
        &=&\pi^{n/4}e^{-i2t-z\bar{z}/2}.\label{eq:conv-bargmann}
\end{eqnarray}
It was shown during the proof of Proposition~\ref{pr:bargmann} that a
convolution with kernel~\eqref{eq:conv-bargmann} defines the usual value
of the Bargmann projector with kernel~\eqref{eq:kr-bargmann}.
\end{proof}
\begin{cor}\label{co:bergman}\textup{\cite[3.1.2]{Rudin80}} The
orthogonal Bergman projector on the space of analytic functions on unit
ball $\Space{B}{}\in \Space{C}{n}$  has
the kernel
\begin{displaymath}
K(\zeta,\upsilon)=(1-\scalar{\zeta}{\upsilon})^{-n-1},
\end{displaymath}
where $\scalar{\zeta}{\upsilon}=\sum_1^n\zeta_j \bar{\upsilon}_j$ is the
scalar product at \Space{C}{n}.
\end{cor}
\begin{proof} We will only rewrite material of Chapters~2 and~3
from~\cite{Rudin80} using our vocabulary. The group of biholomorphic
automorphisms $\object{Aut}(\Space{B}{} )$ of the unit ball \Space{B}{}
acts on \Space{B}{} transitively. For any $\phi\in \object{Aut}
(\Space{B}{})$ there is a unitary operator associated
by~\eqref{eq:unit-trans} and defined by the
formula~\cite[2.2.6(i)]{Rudin80}:
\begin{equation}\label{eq:unit-ball}
[T_\phi f](\zeta)=f(\phi(\zeta))\left(\frac{\sqrt{1-\modulus{\alpha
}^2}}{1-\scalar{\zeta}{\alpha }}\right)^{n+1},
\end{equation}
where $\zeta\in\Space{B}{},\ \alpha =\phi^{-1}(0),\ f(\zeta)\in
\FSpace{L}{2}(\Space{B}{})$. Operator $T_\phi$ from~\eqref{eq:unit-ball}
obviously preserve the space $\FSpace{H}{2}(G)$ of square-integrable
holomorphic functions on \Space{B}{}. The homogeneous space
$\object{Aut}(\Space{B}{})/G_0$ may be identified with
\Space{B}{}~\cite[2.2.5]{Rudin80}. To distinguish points of these two
sets we will denote points of
$B=\object{Aut}(\Space{B}{})/G_0\cong\Space{B}{}$ by Roman letters
(like $a,u,z$) and points of \Space{B}{} itself by Greek letters
($\alpha,
\upsilon, \zeta$ correspondingly). We also always assume, that $a=\alpha,
u=\upsilon, z=\zeta$ under the
mentioned identification.

We select the function $f_0(\zeta)\equiv 1$ as
the vacuum vector.
The mean value formula~\cite[3.1.1(2)]{Rudin80} gives us:
\begin{eqnarray}
\widetilde{f}(a)&=&\scalar{f(\zeta)}{T_\phi f_0} \nonumber \\
      &=&\scalar{T_{\phi^{-1}}f(\zeta)}{ f_0} \nonumber \\
      &=& \int_{\Space{B}{}}  f(\phi(\zeta))\left(\frac{\sqrt{1-
\modulus{\alpha }^2}}{1-\scalar{\zeta}{\alpha }}\right)^{n+1}
\,d\nu(\zeta) \nonumber \\
      &=& f(\phi(0))\left(\frac{\sqrt{1-\modulus{\alpha }^2}}{1-
\scalar{0}{\alpha }}\right)^{n+1} \nonumber \\
      &=& f(a)(\sqrt{1-\modulus{\alpha }^2})^{n+1} ,
\label{eq:f-trans}
\end{eqnarray}
where $a=\alpha =\phi(0),\ \phi\in B$ and
$\widetilde{f}(a)\in\FSpace{L}{2}(B)$. Particulary
\begin{eqnarray}
\widetilde{f}_0(a)&=&(\sqrt{1-\modulus{\alpha  }^2})^{n+1} \nonumber, \\
\widetilde{f}_0(zu)&=&\left(\frac{\sqrt{1-\modulus{\zeta }^2}\sqrt{1-
\modulus{\upsilon }^2}}{1-\scalar{\upsilon}{\zeta}}\right)^{n+1}
\label{eq:f0-trans},
\end{eqnarray}
An invariant measure on $B$ is given~\cite[2.2.6(2)]{Rudin80} by the
expression:
\begin{equation}\label{eq:ball-measure}
d\mu(z)=\frac{d\nu(\zeta)}{(1-\modulus{\zeta}^2)^{n+1}},
\end{equation}
where $d\nu(\zeta)$ is the usual Lebesgue measure on
$B\cong\Space{B}{}$. We
will substitute expressions from~\eqref{eq:f-trans}, \eqref{eq:f0-trans}
and~\eqref{eq:ball-measure} to the reproducing
formula~\eqref{eq:l-reproduce}:
\begin{eqnarray}
\widetilde{f}(u)&=&f(\upsilon)(\sqrt{1-\modulus{\upsilon }^2})^{n+1}
\label{eq:most}\\
     &=&\int_B \widetilde{f}(z) \widetilde{f}_0(z^{-1}u)\, d\mu(z)
\nonumber \\
     &\stackrel{(*)}{=}&\int_B \widetilde{f}(z) \widetilde{f}_0(zu)\,
d\mu(z) \nonumber \\
     &=&\int_{\Space{B}{}}  f(\zeta)(\sqrt{1-\modulus{\zeta }^2})^{n+1}
\left(\frac{\sqrt{1-\modulus{\zeta }^2}\sqrt{1-\modulus{\upsilon
}^2}}{1-\scalar{\upsilon}{\zeta}}\right)^{n+1} \frac{d\nu(\zeta)}{(1-
\modulus{\zeta}^2)^{n+1}} \nonumber \\
     &=&\int_{\Space{B}{}}  f(\zeta)  \left(\frac{\sqrt{1-
\modulus{\upsilon }^2}}{1-\scalar{\upsilon}{\zeta}}\right)^{n+1}\,
d\nu(\zeta) \nonumber \\
     &=&(\sqrt{1-\modulus{\upsilon }^2})^{n+1}\int_{\Space{B}{}}
\frac{f(\zeta)}{(1-\scalar{\upsilon}{\zeta})^{n+1}}\,
d\nu(\zeta) \label{eq:almost}.
\end{eqnarray}
Here transformation $(*)$ is possible because every element of $B$
is an involution~\cite[2.2.2(v)]{Rudin80}.
It immediately follows from the comparison of~\eqref{eq:most}
and~\eqref{eq:almost} that:
\begin{displaymath}
f(\upsilon)=\int_{\Space{B}{}}    \frac{f(\zeta)}{(1-
\scalar{\upsilon}{\zeta})^{n+1}}\, d\nu(\zeta).
\end{displaymath}
The last formula is the integral representation
with the Bergman kernel for holomorphic functions on unit ball in
\Space{C}{n}.
\end{proof}
\begin{cor}\label{co:szego}\textup{\cite{Gindikin64}} The orthogonal
projector Szeg\"o on the boundary
$b\Space{U}{n} $ of the upper half-space in \Space{C}{n+1} has the
kernel
\begin{displaymath}
S(z,w)=(\frac{i}{2}(\bar{w}_{n+1} -z_{n+1})-\sum_{j=1}^nz_j\bar{w}_j)^{-
n-1}.
\end{displaymath}
\end{cor}
\begin{proof} It is well known~\cite{Gindikin64,GreSte77,Stein93} and
was described at Example~\ref{ex:toeplitz}, that there is a unitary
representation of the Heisenberg
group \Heisen{n} as the simply transitive acting group of shift on
$b\Space{U}{n} $ (see~\cite[Chap.~XII, \S~1.4]{Stein93}):
\begin{equation}\label{eq:u-shifts}
(\zeta,t): (z',z_{n+1}) \mapsto (z'+\zeta,
z_{n+1}+t+2i\scalar{z'}{\zeta} +i \modulus{\zeta}^2),
\end{equation}
where $(\zeta,t)\in\Heisen{n},  z=(z',z_{n+1})\in\Space{C}{n+1},
\zeta,z'\in\Space{C}{n}, t\in \Space{R}{}$. We again apply the general
scheme from Example~\ref{ex:transform}. This gives an identification of
\Heisen{n} and $b\Space{U}{n}$ and \Heisen{n} act on
$b\Space{U}{n}\cong\Heisen{n}$ by left group shifts. Left invariant
vector fields are exactly the {\em tangential Cauchy-Riemann
equations\/} for holomorphic functions on \Space{U}{n}.
Shifts~\eqref{eq:u-shifts} commute with the tangential Cauchy-Riemann
equations and thus preserve the {Hardy space}
$\FSpace{H}{2}(b\Space{U}{n} ) $ of boundary values of functions holomorphic on
\Space{U}{n}.

As vacuum vector we select the function $f_0(z)=(iz_{n+1})^{-n-
1}\in\FSpace{H}{2}(b\Space{U}{n} )$. Then the Szeg\"o projector
$P:\FSpace{l}{2}(b\Space{U}{n} ) \rightarrow \FSpace{H}{2}(b\Space{U}{n}
) $ is the right convolution on $\Heisen{n} \cong b\Space{U}{n}$ with
$f_0(z)$ and thus should have the kernel (see the group low
formula~\eqref{eq:h-low} for \Heisen{n})
\begin{displaymath}
S(z,w)=(\frac{i}{2}(\bar{w}_{n+1} -z_{n+1})-\sum_{j=1}^nz_j\bar{w}_j)^{-
n-1}.
\end{displaymath}

Reader may ask, {\em why have we selected such a vacuum vector?\/} The
answer is: for a {\em simplicity\/} reason. Indeed, the {\em Cayley
transform\/} (\cite[Chap.~XII, \S~1.2]{Stein93}
and~\cite[\S~2.3]{Rudin80})
\begin{equation}\label{eq:cayley}
C(z)=i\frac{e_{n+1}+z}{1-z_{n+1}},\
e_{n+1}=(0,\ldots,0,1)\in\Space{C}{n+1}
\end{equation}
establishes a biholomorphic map from unit ball $\Space{B}{}\in
\Space{C}{n+1}$ to domain \Space{U}{n}. We can construct an isometrical
isomorphism of the Hilbert spaces
$\FSpace{H}{2}(\Space{S}{2n+1})$ and $\FSpace{H}{2}(\Space{U}{n})$ based
on~\eqref{eq:cayley}
\begin{equation}\label{eq:c-iso}
f(z) \mapsto [Cf](z)= f(C(z))\frac{-2i^{n+1}z_{n+1}}{(1-
z_{n+1})^{n+2}},\ f\in \FSpace{H}{2}(\Space{U}{n}),\
[Cf]\in\FSpace{H}{2}(\Space{S}{2n+1}).
\end{equation}
Then the vacuum vector $f_0=(iz_{n+1})^{-n-1}$ is the image of function
$\widetilde{f}_0(w)=(-2i/(w-i))^{n+2}\in\FSpace{H}{2}(\Space{S}{2n+1})$
under transformation~\eqref{eq:c-iso}. It seems to be one from the
simplest functions from $\FSpace{H}{2}(\Space{S}{2n+1})$ with
singularities on $\Space{S}{2n+1}$.
\end{proof}

\subsection{Connections with Relative Convolutions}\label{ss:connect}
Now we return to relative convolutions and will show their connections
with coherent states. For any operator
$A:\FSpace{L}{2}(G,d\mu)\rightarrow\FSpace{L}{2}(G,d\mu)$ we can
construct the {\em Toeplitz-like\/} operator
$P_A=PA:\FSpace{H}{2}(G,d\mu)\rightarrow\FSpace{H}{2}(G,d\mu)$. Of
course, using the isomorphism $H\cong\FSpace{H}{2}(G,d\mu)$ we can think
about $P_A$ as an operator $P_A:H\rightarrow H$. Particulary, operators
of group convolution on $G$ will induce relative convolutions on $H$
(see Lemma~\ref{le:category}). Thus we again have a direct way for
applications of harmonic analysis in every problem concerning coherent
states. Although, coherent states are very useful in physics we will
stop here\footnote{However, let us remind again that
Example~\ref{ex:fock} forms an interesting
application to physics.} and will only develop this theme in the next
Example connected with wavelets.

\section{Applications to Physics and Signal Theory}\label{se:physics}
We are going to consider some Examples connected with physics, but our
division between mathematics and physics is so fragile as in the real life.
\begin{example}\label{ex:wavelets} Let us consider the ``$ax+b$
group''~\cite[\S~7.1]{MTaylor86} of affine transformations of the real
line.  We will denote this group by $A$ and its Lie algebra by
\algebra{a}. We will consider their operation on the real line
$S=\Space{R}{} $. \algebra{a} is spanned by two vector fields
$X_s=\frac{1}{i}\frac{\partial }{\partial y}$ (which generate shifts)
and $X_d=y\frac{1}{i}\frac{\partial }{\partial y}$ (which generate
dilations). Their commutators are $[X_s,X_d]=X_s$. Then
transformation~\eqref{eq:embedding} takes the form
\begin{eqnarray}
\widetilde{f}(x_1,x_2)&=&\frac{1}{2\pi}\int_{\Space{R}{}} \bar{f}(y)\,
e^{i(x_1 X_s+ x_2 X_d)}f_0(y)\,dy\nonumber \\
   &=&\frac{1}{2\pi}\int_{\Space{R}{}} \bar{f}(y)\,
e^{-x_2/2}f_0(e^{-x^2}y-x_1)\,dy. \label{eq:wavelet}
\end{eqnarray}
The last line is easily recognized as the {\em wavelet
transform\/}~\cite{HeilWaln89}.

Similar expression in the spirit of
Definition~\ref{de:coherent1} for the {\em Gabor
transform\/}~\cite{HeilWaln89} may be
obtained if we replace $A$ by the Heisenberg group \Heisen{3}. Then, as
was shown early, different signal alterations constructed in signal
theory are relative convolution on $S=\Space{R}{}$ induced  by
\Heisen{3} or the meta-Heisenberg group from~\cite{Folland94}. For
example,
signal filtration may be presented at the Gabor representation as a
multiplication by characteristic functions of the wished time and frequency
intervals.

Using relative convolutions we can coherently introduce wavelets-like
transform for every semi-direct product~\cite[\S~5.3]{MTaylor86}  of Lie
group and Abelian one~\cite{BernTayl94}. In view of applications to the
signal theory it
seems interesting to start from the Heisenberg group and its dilations.
The one-parameter group $D=\{\delta_\tau \such \tau\in\Space{R}{} \}$ of
dilations of the Heisenberg group is given by the formula
\begin{equation}\label{eq:dilation}
\delta_\tau(t,z)=(e^{2\tau}t,e^\tau z),\ (t,z)\in\Heisen{n} \cong
\Space{R}{}\times\Space{R}{2n}
\end{equation}
and has the one-dimensional Lie algebra spanned on the vector field
\begin{equation}\label{eq:frame-dilat}
X_h=\frac{1}{i}(2t\frac{\partial }{\partial t} + z\frac{\partial
}{\partial z}).
\end{equation}
If we introduce now the relative convolutions for Lie algebra generated
by vector fields $X^l_j$ from~\eqref{eq:frame-left} and $X_h$
from~\eqref{eq:frame-dilat} then we obtain the {\em Heisenberg
Gabor-like transform\/}, which should be useful to analyze of the
radar ambiguity function~\cite[\S~1.4]{Folland89}.
\end{example}
\begin{example}\label{ex:quantum} We are going to describe {\em
group quantization\/} from paper~\cite{Kisil94d}.  The usual
``quantization'' means
some (more or less complete) set of rules for the construction of a
quantum algebra from the classical description of a physical system. The
group quantization is based on the Hamilton description and consists of
the following steps.
\begin{enumerate}

\item Let $\Omega=\{x_j\}, 1\leq j\leq N$ be a set of physical
quantities defining state of a classical system. Observables are real
valued functions on the states.

The most known and important case is the set $\{x_j=q_j,x_{j+n}=p_j\},\
1\leq j\leq n, N=2n$ of coordinates and impulses of classical particle
with $n$ degrees of freedom. Observables are real valued functions on
\Space{R}{2n}. This example will be our main illustration during the
present consideration.

\item We will complete the set $\Omega$ till $\bar{\Omega}$ by
additional quantities ${x_j}, N<j\leq \bar{N}$, such that $\bar{\Omega}$
will form the smallest algebra containing $\Omega$ and closed under the
Poisson bracket:
\begin{displaymath}
\{x_i,x_j\}\in \bar{\Omega},\ \mbox{ for all } x_i,x_j\in \bar{\Omega}.
\end{displaymath}

In the case of a particle we should add the function $x_{2n+1}=1$, which
is equal to unit identically and one obtains the famous relations
($\bar{N}=2n+1$)
\begin{equation}\label{eq:poisson}
\{x_j,x_{j+n}\}=-\{x_{j+n},x_j\}=x_{2n+1}
\end{equation}
and all other Poisson brackets are equal to zero.

\item We form a $\bar{N}$-dimension Lie algebra $\algebra{p}$ with a
frame $\{\widehat{x}_j\},\ 1\leq j \leq \bar{N}$ with the formal mapping
$\hat{}: x_j\mapsto\widehat{x}_j$. Commutators of frame vectors of
\algebra{p} are
formally
defined throughout the formula
\begin{equation}\label{eq:hat-poisson}
[\widehat{x}_i,\widehat{x}_j]=\widehat{\{x_i,x_j\}}
\end{equation}
and we extend the commutator on whole algebra by the linearity.

For a particle this step give us the Lie algebra
$\algebra{h}_n$ of the Heisenberg group (compare~\eqref{eq:heisen-comm}
and~\eqref{eq:poisson}).

\item We introduce an algebra $\algebra{P}$ of relative
convolutions~\eqref{eq:relative} induced by \algebra{p}. These
operators are {\em observables\/} in the group quantization and by
analogy to classic case they may be treated as functions of
$\widehat{x}_j$ (see Remark~\ref{re:function}). A set $S$
which algebra \algebra{p} acts on and type of kernels are depending on
physically determining constraints. The family of all one-dimensional
representations of \algebra{P} is called {\em classical\/} mechanics and
different noncommutative representations correspond to {\em quantum\/}
descriptions with the different {\em Planck constants\/}.

For particle we have the following opportunities:
\begin{enumerate}
\item $S=\Space{R}{n},\ \widehat{x}_j=X_j=M_{q_j},\
\widehat{x}_{j+n}=\hbar\frac{1}{i}\frac{\partial }{\partial q_j}$,
relative convolutions are PDO from Example~\ref{ex:pdo} and we have
obtained the {\em Dirac--Heisenberg--Schr\"odinger--Weyl quantization}
by PDO.

\item $S=\Space{R}{2n},\ \widehat{x}_j=X_j=M_{q_j},\
\widehat{x}_{j+n}=M_{p_j}$, relative convolutions are operators of
multiplication by functions (or just functions) from
Example~\ref{ex:multiplication}  and we have obtained the usual
classical description, which we have started from.

\item $S=\Heisen{n},\ \widehat{x}_j=X^{l(r)}_j,\ 0\leq j\leq 2n+1$ and
relative convolutions form the group convolution algebra on \Heisen{n}.
This description (so-called {\em plain mechanics\/}) contain both the
{\em quantum\/} and {\em classical\/} ones with the natural realization
of the {\em correspondence principle\/} (see~\cite{Kisil94d} for
details).
\end{enumerate}
\end{enumerate}

Group quantization is straightforward enough and obviously
preserves the symmetry group of the classical system under
investigation. Moreover, there is also other advantages of
the proposed quantization, which distinguish it from
the already known ones.
\begin{itemize}
\item In contrary to the {\em operator quantization\/} of
Berezin~\cite{Berezin74} and the {\em geometrical quantization\/} of
Kirillov--Souriau--Konstant~\cite{Woodhouse80} we should not introduce
{\em a~priori\/} any Planck constants. Moreover, during the posterior
analysis of relative convolution algebra representations a parameter
corresponding to the Planck constant will appear naturally. By the way,
a {\em set of Planck constants\/} should not necessary belong to
$[0,+\infty[$ and may form more complicated topological spaces.

\item The problem of ordering of noncommutative quantities
$\widehat{q}_j$ and $\widehat{p}_j$ does not occur under group
quantization. Correspondence
\begin{displaymath}
\mbox{function } k(x)\ \rightarrow \mbox{ convolution
with the kernel }k(x)
\end{displaymath}
 is direct enough even for noncommutative groups.
Meanwhile, in other quantization the ``painful question of
ordering''~\cite{Berezin71} has generated many different answers:
the $\widehat{q}\widehat{p}$-quantization, the
$\widehat{p}\widehat{q}$-quantization, the Weyl-symmetrical,
the Wick and  the anti-Wick
(Berezin) quantization.
\end{itemize}
\end{example}
\begin{rem} The presented group quantization has deep roots in the
quantization procedure of Dirac~\cite{Dirac67}. The main differences are
\begin{itemize}
\item We recognize the Heisenberg commutation
relations~\eqref{eq:heisen-comm} only as a particular
case among other possibilities. However, due to Theorem~\ref{th:pdo} they
play the fundamental role.
\item We do not look only for irreducible representations of commutation
relations.
\end{itemize}
\end{rem}

\section{Conclusion}
The paper tried to illustrate {\em how the systematical usage of
harmonic analysis in various applications may be useful for both:
the analysis and applications\/}. It seems, that relative convolutions form
an appropriate tool for this purpose.

Given Examples from different fields of mathematics and physics made
reasonable studying of relative convolutions. Moreover, we  have
repeatedly met nilpotent Lie groups (and particularly the Heisenberg)
group within important applications, so our primary interest in such
groups should be excused.

\bibliographystyle{amsplain}
\bibliography{mrabbrev,akisil,analyse,aphysics}

\end{document}